\documentclass[]{amsart}

\def\bi{\mathbb{I}}
\def\bj{\mathbb{J}}
\def\bk{\mathbb{K}}
\def\bn{\mathbb{N}}
\def\br{\mathbb{R}}
\def\bc{\mathbb{C}}

\def\h{\mathcal{H}}
\def\k{\mathcal{K}}
\def\l{\mathcal{L}}
\def\g{\mathcal{G}}
\DeclareMathOperator{\B}{\rm B}

\newcommand{\nor}[1]{#1_{\rm n}}
\newcommand{\adb}[1]{#1^{_{A}\sharp_{B}}}
\newcommand{\md}[1]{#1^{\natural}}
\newcommand{\pmd}[1]{#1^{\natural_{\rm p}}}
\newcommand{\pmp}[1]{#1_{\natural_{\rm p}}}

\newcommand{\mbd}[1]{#1^{\natural\natural}}
\newcommand{\pmbd}[1]{#1^{\natural_{\rm p}\natural_{\rm p}}}

\newcommand{\mpd}[1]{#1_{\natural}}

\newcommand{\mdbxy}{(X\stackrel{h}{\otimes}_BY)^{\natural_{B{\rm nor}}}}
\newcommand{\oxv}{\omega_{x,v}}
\newcommand{\rxv}{\rho_{x,v}}
\newcommand{\toxv}{\tilde{\omega}_{x,v}}
\newcommand{\trxv}{\tilde{\rho}_{x,v}}

\newcommand{\du}[1]{#1^{\sharp}}
\newcommand{\pd}[1]{#1_{\sharp}}
\newcommand{\bdu}[1]{#1^{\sharp\sharp}}
\newcommand{\com}[1]{#1^{\prime}}
\newcommand{\ort}[1]{#1^{\perp}}

\newcommand{\norm}[1]{\Vert#1\Vert}
\newcommand{\cbnorm}[1]{\Vert#1\Vert_{{\rm cb}}}
\newcommand{\inner}[2]{\langle #1,#2\rangle}
\newcommand{\sumjn}{\sum_{j=1}^n}

\newcommand{\sumj}{\sum_{j\in\mathbb{J}}}

\newcommand{\bb}[2]{{\rm B}(#1,#2)}
\newcommand{\abbb}[2]{{\rm B}_A(#1,#2)_B}
\newcommand{\anbbb}[2]{{\rm N}_A(#1,#2)_B}

\newcommand{\bh}{{\rm B}(\mathcal{H})}
\def\bk{{\rm B}(\mathcal{K})}
\newcommand{\bl}{{\rm B}(\mathcal{L})}

\newcommand{\Ll}{{\rm L}}
\newcommand{\ta}{\tilde{A}}
\newcommand{\tb}{\tilde{B}}

\newcommand{\bkh}{{\rm B}(\mathcal{K},\mathcal{H})}

\newcommand{\trm}{{\rm T}_m}
\newcommand{\cbb}[2]{{\rm CB}(#1,#2)}

\newcommand{\acbbb}[2]{{\rm CB}_A(#1,#2)_B}

\newcommand{\ncbb}[2]{{\rm NCB}(#1,#2)}
\newcommand{\ancbbb}[2]{{\rm NCB}_{A}(#1,#2)_{B}}

\newcommand{\acbb}[1]{{\rm CB}_A(#1)_B}
\newcommand{\ancbb}[1]{{\rm NCB}_A(#1)_B}

\newcommand{\matn}[1]{{\rm M}_{n}(#1)}
\newcommand{\matm}[1]{{\rm M}_{m}(#1)}
\newcommand{\matj}[1]{{\rm M}_{\mathbb{J}}(#1)}

\newcommand{\rown}[1]{{\rm R}_{n}(#1)}
\newcommand{\coln}[1]{{\rm C}_{n}(#1)}
\newcommand{\rowj}[1]{{\rm R}_{\mathbb{J}}(#1)}
\newcommand{\colj}[1]{{\rm C}_{\mathbb{J}}(#1)}

\newcommand{\rowi}[1]{{\rm R}_{\mathbb{I}}(#1)}
\newcommand{\coli}[1]{{\rm C}_{\mathbb{I}}(#1)}
\newcommand{\mat}[3]{{\rm M}_{#1,#2}(#3)}

\newcommand{\weakc}[1]{\overline{#1}}
\newcommand{\hg}{\stackrel{h}{\otimes}}
\newcommand{\ehg}{\stackrel{eh}{\otimes}}

\newcommand{\spac}{\check{\otimes}}
\newcommand{\nhg}{\stackrel{\sigma h}{\otimes}}
\newcommand{\opr}{\hat{\otimes}}

\newcommand{\Max}[1]{{\rm MAX}(#1)}

\newcommand{\maxab}[1]{{\rm MAX}_A(#1)_B}
\newcommand{\maxn}[1]{{\rm MAXN}_A(#1)_B}
\newcommand{\maxnd}[1]{{\rm MAXND}_A(#1)_B}
\newcommand{\minab}[1]{{\rm MIN}_A(#1)_B}

\newcommand{\m}[1]{\|#1\|_{_A{\rm m}_B}}
\newcommand{\M}[1]{\|#1\|_{_A{\rm M}_B}}
\newcommand{\Mn}[1]{\|#1\|_{_A{\rm MN}_B}}
\newcommand{\Mnd}[1]{\|#1\|_{_A{\rm MND}_B}}
\newcommand{\os}{{\rm OS}}
\newcommand{\abb}{{_A{\rm BM}_B}}
\newcommand{\adbb}{{_A{\rm DBM}_B}}
\newcommand{\aob}{{_A{\rm OM}_B}}
\newcommand{\anob}{{_A{\rm NOM}_B}}
\newcommand{\adob}{{_A{\rm DOM}_B}}
\newcommand{\andob}{{_A{\rm NDOM}_B}}
\newcommand{\arb}{{_A{\rm RM}_B}}
\newcommand{\anrb}{{_A{\rm NRM}_B}}
\newcommand{\andrb}{{_A{\rm NDRM}_B}}
\newcommand{\asob}{{_A{\rm SOM}_B}}

\newcommand{\coc}{{{\rm COM}_C}}

\newcommand{\cnoc}{{{\rm CNOM}_C}}

\DeclareMathOperator{\essup}{\rm essup}
\DeclareMathOperator{\eop}{{\rm ess}\oplus_{t\in\Delta}}

\newtheorem{theorem}{Theorem}[section]
\newtheorem{lemma}[theorem]{Lemma}
\newtheorem{corollary}[theorem]{Corollary}
\newtheorem{proposition}[theorem]{Proposition}
\theoremstyle{remark}
\newtheorem{remark}[theorem]{Remark}
\theoremstyle{definition}
\newtheorem{definition}[theorem]{Definition}
\newtheorem{example}[theorem]{Example}
\numberwithin{equation}{section}

\begin{document}

\title[]{Duality and normal parts of operator modules} 
\author{Bojan Magajna} 
\address{Department of Mathematics\\ University of Ljubljana\\
Jadranska 19\\ Ljubljana 1000\\ Slovenia}
\email{Bojan.Magajna@fmf.uni-lj.si}

\thanks{Supported by the Ministry of Science and Education of Slovenia.}
\keywords{Operator bimodule, von Neumann algebra, 
relative tensor products.}

\subjclass[2000]{Primary 46L07, 46H25; Secondary 47L25}

\begin{abstract} For an operator bimodule $X$ over von Neumann algebras 
$A\subseteq\bh$ and $B\subseteq\bk$, 
the space of all completely bounded $A,B$-bimodule maps
from $X$ into $\bkh$, is the bimodule dual of $X$.
Basic duality theory is developed with a particular attention to 
the Haagerup tensor product over von Neumann algebras. To $X$  
a normal operator bimodule $\nor{X}$ is associated so that completely bounded 
$A,B$-bimodule maps from $X$ into  normal operator 
bimodules factorize uniquely through $\nor{X}$. A construction of $\nor{X}$ in
terms of biduals of $X$, $A$ and $B$ is presented.
Various operator bimodule structures are considered on a 
Banach bimodule admitting a normal such structure.  

\end{abstract}

\maketitle

\section{Introduction}

One of the aims of this article is to show that the  classical duality 
theory of Banach
spaces and the more recent duality of operator spaces (\cite{BLM}, \cite{ER},
\cite{P}, \cite{Pi}) effectively extends
to the situation where a Banach space is replaced by a normal operator bimodule $X$
over von Neumann algebras $A$ and $B$.  The role of
the dual is played by the $\com{A},\com{B}$-bimodule $\md{X}$ consisting of all completely bounded
$A,B$-bimodule maps from $X$ into $\bkh$, where $\h$ and $\k$ are proper Hilbert
modules over $A$ and $B$, respectively.  Among the basic tools (or motivations)
for such an extension are the operator valued Hahn - Banach and bipolar theorems 
(\cite{Ar}, \cite{EW}, \cite{W1}). Some general aspects of duality of operator bimodules were considered
also in \cite{N} and \cite{Po}, \cite{AP}. Here we study mainly bimodules over von 
Neumann algebras and emphasize the normality considerations. We shall explain briefly
an application to  W$^*$-correspondences (Section 3).

In Section 2 we collect definitions of various (known) classes of bimodules,
introduce  abbreviations for their names and summarize some  preliminary 
results.  

In Section 3 we develop our basic technique and prove some typical duality 
theorems. In the formulation of results we 
are guided by classical functional analysis, but since the range  of 
`functionals' here is $\bkh$ instead of $\bc$, the proofs of 
main results require methods 
completely different from the classical ones. Our starting point will be a duality
result for the Haagerup tensor product of normal operator bimodules
(Theorem \ref{t32}), which extends the duality for the usual Haagerup tensor product
of operator spaces obtained by Blecher and Smith \cite{BS}. This will enable us to
relate the bimodule duals to the usual operator space duals. Many classical 
results can be generalized at least to strong bimodules.
In Section 3 we also consider very briefly relations between the properties of a 
given bimodule
map $T$ and its bimodule adjoint $\md{T}$.

Because of the central role
of the extended module Haagerup tensor product,  we relate in Section 4 its 
bimodule dual to the normal version of the Haagerup tensor product studied 
by Effros and Ruan in \cite{ER3}. We also describe the module versions of
the extended and the normal Haagerup tensor products of two von Neumann algebras
over a common von Neumann subalgebra as concrete spaces of operators and thus
generalize some results of Blecher - Smith \cite{BS}  and  Effros - Kishimoto 
\cite{EK}. 
 
For a general operator $A,B$-bimodule $X$ we shall show that the
closure $\nor{X}$ of the image of $X$ in its bimodule bidual
$\mbd{X}$ is a normal operator $A,B$-bimodule having
the following universal
property: for each completely bounded $A,B$-bimodule map
$\phi$ from $X$ to a normal operator $A,B$-bimodule $Y$ there exists a unique
$A,B$-bimodule map $\tilde{\phi}$ from $\nor{X}$ into $Y$ such that
$\phi=\tilde{\phi}\iota$, where $\iota$ is the canonical map from
$X$ into $\nor{X}$. The bimodule $\nor{X}$ is described in Section 5 in
an alternative way and is called the normal part of $X$,
although in general it is not contained in $X$. We also consider how
$\mbd{X}$ sits completely isometrically in the operator space
bidual $\bdu{X}$ of $X$ (Theorem \ref{t59}),
where $\bdu{X}$ is equipped with the canonical
normal operator bimodule structure over the universal von Neumann
envelopes of $A$ and $B$.

In Section 6 the discussion is specialized to central bimodules over an abelian von Neumann algebra
$C$.
A $C$-bimodule $X$ is called central if $xc=cx$ 
for all $c\in C$ and $x\in X$. First we observe that a central operator $C$-bimodule $X$ is 
normal if and only if for each $x\in\matn{X}$ the function
$\Delta\ni t \mapsto\norm{x(t)}$ on the spectrum $\Delta$ of $C$ is continuous, where $x(t)$ is the coset of
$x$ in $\matn{X}/[({\rm ker\, t})\matn{X}]$. Then we characterize concretely the normal part of a central
$C$-bimodule $X$ in terms of its decomposition along $\Delta$ (Theorem \ref{t511}). 
We also prove that 
for a strong central $C$-bimodule $X$ 
and a subbimodule $Y$ in $X$ the quotient $X/Y$ is
normal if and only if $Y$ is strong.

If a Banach $A,B$-bimodule $X$ over von Neumann algebras $A$ and $B$ admits a norm structure of 
a {\em normal} operator $A,B$-bimodule then  the maximal  such structure
$\maxn{X}$ turns out to be different from the maximal operator $A,B$-bimodule 
$\maxab{X}$ on $X$ (Section 7).  Indeed, $\maxn{X}$ is  just the normal part
of $\maxab{X}$. Even more surprisingly,  if $X$ admits
a structure of a normal {\em dual} operator
$A,B$-bimodule, the maximal such structure, $\maxnd{X}$, is different from
$\maxn{X}$.  This  provides new examples of operator spaces
which are duals as Banach spaces, say $\du{V}$, but without any  
operator space predual  on $V$. 
Earlier such examples are in \cite{Le} and \cite{EOR}.

\section{Basic classes of bimodules, notation and other preliminaries}

Throughout the paper
$A,B$ and $C$ will be C$^*$-algebras with unit $1$, in fact von Neumann algebras 
most of the time. By a 
Banach $A,B$-bimodule we mean a Banach space
$X$ which is an $A,B$-bimodule such that $1x=x=x1$ and $\norm{axb}\leq\norm{a}
\norm{x}\norm{b}$ for all $a\in A$, $b\in B$ and $x\in X$. The class of 
all such bimodules is denoted by $\abb$,  and the space of all
bounded $A,B$-bimodule maps from $X$ to $Y$ by 
$\abbb{X}{Y}$. 

A Hilbert $A$-module is just a Hilbert space $\h$ 
together with a $*$-representation
$\pi$ of $A$ on $\h$. We shall always assume that $\pi(1)=1_{\h}$. If $A$ is a von Neumann algebra and $\pi$ is normal
then $\h$ is called {\em normal}. 
If $\pi$ is injective, $\h$ is called {\em faithful}.
If $\pi$ is cyclic then $\h$ is called 
{\em cyclic}.  If each finite subset of $\h$
is contained in a closed cyclic submodule $[A\xi]$ then
$\h$ is  {\em locally cyclic}. The importance of such modules originates from
a well known result of \cite{S} recalled in Theorem \ref{acb}. Over
a von Neumann algebra $A$ there exists a canonical 
Hilbert module $\h$, called the standard form. We shall only need 
a property of this module abstracted in the following 

\begin{definition} {\em A proper module} over a von Neumann algebra $A$ is
a faithful normal Hilbert $A$-module $\h$ such that all normal states on $A$
and on $\com{A}$ (the commutant of $A$ in $\bh$) are vector states (that is,
of the form $x\mapsto\inner{x\xi}{\xi}$, $\xi\in\h$). 
\end{definition}

Note that a proper $A$-module $\h$ contains (up to a unitary equivalence) all normal 
cyclic Hilbert $A$-modules (since all normal states on $A$ come from vectors in $\h$)
and is locally cyclic by \cite[2.3]{S}. 
Since for a separable $\h$ locally cyclic vectors are cyclic by \cite[2.7]{HW},
it follows from \cite[7.2.9]{KR} that  a proper separable 
$A$-module is essentially just the standard form.

For operator spaces $X$ and $Y$,  $\cbb{X}{Y}$ denotes the set of completely bounded
linear maps from $X$ to $Y$. Occasionally we shall use the notation $\os$ for
the class of operator spaces. If $X,Y\in\os\cap\abb$, let
$\acbbb{X}{Y}=\abbb{X}{Y}\cap\cbb{X}{Y}$. We are now going to recall the definitions
of various classes of operator modules. We will follow the usual terminology, but
since some classes of modules have very long names (such as `normal dual operator
$A,B$-bimodules') and appear repeatedly, it will be convenient to introduce notation for 
them. 

\begin{definition}\label{no}(i) The class $\aob$ of {\em operator $A,B$-bimodules} 
consists of all bimodules $X\in\abb\cap\os$ which can be completely isometrically
and homomorphically represented in a $\bh$. In other words, for some Hilbert 
module $\h$ over $A$ and $B$ the space $\acbbb{X}{\bh}$ contains a complete isometry.

(ii) If in (i) $A$ and $B$ are von Neumann algebras and  $\h$ can be chosen to be
normal over $A$ and $B$, then  $X$ is a {\em normal operator $A,B$-bimodule} 
($X\in\anob$).
\end{definition}

As we shall observe below, the class of normal operator modules contains
normal Hilbert modules and will play an important role here.

Operator bimodules are characterized by the CES 
theorem \cite{CES}, which was later generalized and sharpened  (\cite[Section 5]{B4},  \cite[Section 4.6]{BLM},
\cite[Chapter 16]{P}).
\begin{theorem}\label{om}\cite{BLM}, \cite{P} A bimodule $X\in\os\cap\abb$ is an operator bimodule
if and only if $\matn{X}$ 
is a 
Banach $\matn{A},\matn{B}$-bimodule for each $n=1,2\ldots$. 

Given $X\in\os$, there
exist C$^*$-algebras $A_{l}(X)$ and $A_r(X)$ such that $X$ is an operator
$A_l(X),A_r(X)$-bimodule and every left (right) operator $A$-module structure
on $X$ is given by a $*$-homomorphism from $A$ into $A_l(X)$ (into $A_r(X)$).
In particular, if $X$ is a left operator $A$ module and a right operator
$B$-module, then $X$ is automatically an operator $A,B$-bimodule.
\end{theorem}

Normal operator bimodules are characterized by the following result, 
the first part of which is not hard to
deduce also from  Theorem \ref{t59} below (and its proof).
\begin{theorem}\label{tno}\cite[3.3]{M6}, \cite[6.1]{M8} A bimodule 
$X\in\aob$ is normal if and only if for each 
$n\in\bn$ and  $x\in\matn{X}$ the mappings
$\matn{A}\ni a\mapsto\norm{ax}$ and $\matn{B}\ni b\mapsto\norm{xb}$
are weak* lower semicontinuous. If $A$ and $B$ are $\sigma$-finite,
this is the case if and only if for all $x\in\matn{X}$ and
sequences of projections $(e_j)$ and $(f_j)$ increasing to $1$ in $\matn{A}$ 
and $\matn{B}$ (resp.) we have
$\lim_j\norm{e_jx}=\norm{x}=\lim_j\norm{xf_j}$.
\end{theorem}

We recall that a von Neumann algebra $A$ is  $\sigma$-finite if each orthogonal
family of nonzero projections in $A$ is countable.

\begin{definition}\label{db} (i) A {\em dual Banach $A,B$-bimodule} is a dual Banach 
space $X=\du{V}\in\abb$ such that the maps
$X\ni x\mapsto ax$ and $X\ni x\mapsto xb$
are weak* continuous for all $a\in A$ and $b\in B$. Then the
preadjoints of these maps  define a $B,A$-bimodule structure
on $V$. Conversely, for every  $V\in{_B{\rm BM}_A}$,
$X=\du{V}$ becomes a {\em  dual Banach $A,B$-bimodule}
by
$$\inner{axb}{v}=\inner{x}{bva}\ (x\in X,\ v\in V).$$
The category of such bimodules is denoted by $\adbb$
and the space of all weak* continuous (hence bounded) $A,B$-bimodule 
maps  from $X$ to $Y$ by $\anbbb{X}{Y}$.

(ii) If  $A$ and $B$ are von Neumann algebras, a bimodule $X=\du{V}\in\adbb$ is a  
{\em normal dual Banach bimodule} ($X\in{_A{\rm NDBM}_B}$) 
if the maps $A\ni a\mapsto\inner{ax}{v}$ and
$B\ni b\mapsto\inner{xb}{v}$ are weak* continuous for all $x\in X$ and $v\in V$.
\end{definition}

\begin{definition}\label{do}  An operator bimodule $X$ is a
{\em  dual operator $A,B$-bimodule} ($X\in\adob$) if $X$ is the operator space dual 
of some $V\in{_B{\rm BM}_A}\cap\os$ equipped with the $A,B$-bimodule action as in
Definition \ref{db}. 
For such bimodules $X$, $Y$ let $\ancbbb{X}{Y}=\anbbb{X}{Y}\cap\cbb{X}{Y}$. 
\end{definition}

\begin{remark}\label{rh} A Hilbert $A$-module $\h$ is regarded as an operator 
$A$-module by considering $\h$ as the column operator space (\cite{BLM} or
\cite{Pi}). Then $\h$ is dual to the
conjugate Hilbert space $\h^*$ with the row operator space structure and
the right module action 
\begin{equation}\label{hil} \xi^*a=(a^*\xi)^*,\ \ (\xi\in\h,\ a\in A),
\end{equation}
where $\xi^*$ denotes $\xi$ regarded as an element of $\h^*$.
{\em In this paper $\h$ will always mean a column Hilbert space and $\h^*$ the
corresponding operator space dual.}
\end{remark}

\begin{definition}\label{dno} For von Neumann algebras $A$ and $B$, a bimodule
$X\in\adob$
is {\em a normal dual operator $A,B$-bimodule} ($\andob$) if 
there exist a normal Hilbert module $\h$ over $A, B$ and
a complete isometry in $\ancbbb{X}{\bh}$. 
\end{definition}

The original characterization of normal dual operator bimodules \cite{ER1} was 
greatly improved
in \cite[4.1, 4.2]{B3} and \cite[5.4, 5.7]{BEZ}, from which we shall need the following: 

\begin{theorem}\label{thbez}\cite{B3}, \cite{BEZ}   If $X$ is a dual operator space, 
then the algebras
$A_l(X)$ and $A_r(X)$ in Theorem \ref{om} are von Neumann algebras and $X$ is a normal dual operator 
$A_l(X),A_r(X)$-bimodule.  Thus, if $A$ and $B$ are von Neumann algebras and $X\in\adob$, the 
maps $X\ni x\mapsto ax$ ($a\in A$) and $x\mapsto xb$ ($b\in B$) are 
automatically weak* continuous. If 
the maps $A\ni a\mapsto ax$ and $B\ni b\mapsto xb$ are also weak* continuous for each $x\in X$, 
then $X\in\andob$.
\end{theorem}
\begin{remark}\label{rnd} If a dual operator $A,B$-bimodule $X$ satisfies the norm semicontinuity
condition of Theorem \ref{tno} then $X$ is a normal dual operator bimodule (see
\cite[p. 199--200]{M6}), that is $\adob\cap\anob=\andob$.
\end{remark}
For an index set $\bj$ and an $X\in\os$, let $\rowj{X}$ be the set $\mat{1}{\bj}{X}$ of all 
$1\times \bj$ bounded matrices
with the entries in $X$. Similarly, $\colj{X}:=\mat{\bj}{1}{X}$. 
(An $\bi\times\bj$ matrix is 
bounded if the supremum of the norms of its finite submatrices is finite.)

\begin{definition}\label{dstr} A bimodule $X\in\anob$ is called {\em strong} ($X\in\asob$) if  
\begin{equation}\label{str}[a_i][x_{ij}](b_j)=\sum_{i\in\bi,j\in\bj}a_ix_{ij}b_j\in X
\end{equation}
for all  $[a_i]\in\rowi{A}$,
$[x_{ij}]\in\mat{\bi}{\bj}{X}$, $(b_j)\in\colj{A}$ 
and all index sets $\bi$ and $\bj$.  
\end{definition}
As shown in \cite{M5}, it suffices to require the condition (\ref{str}) for
orthogonal families of projections $(a_i)\subseteq A$ and 
$(b_j)\subseteq B$.
Strong bimodules in $\bh$ are characterized
as closed in the $A,B$-topology \cite{M7}, the definition of 
which will not be needed here. For our purposes
it will suffices to note that a  functional $\rho$ on $\bh$  
is $A,B$-continuous if and only 
if $\rho\in\adb{\bh}$, where $\adb{\bh}$ is defined as follows.  

\begin{definition}\label{adb} If $A$ and $B$ are von Neumann algebras and
$X\in\abb$,  let $\adb{X}$ be the subspace of the dual $\du{X}$ of $X$, consisting
of all $\rho\in\du{X}$ such that for each $x\in X$ 
the maps $A\ni a\mapsto\rho(ax)$ and $B\ni b\mapsto\rho(xb)$
are weak* continuous.
\end{definition}

The argument from \cite[4.6]{M5} shows
that bounded bimodule homomorphisms are continuous in the $A,B$-topology. 

Occasionally we shall need a version of the operator bipolar theorem.
A subset $K$ of a bimodule $X\in\abb$ is called {\em  $A,B$-absolutely convex}
if
$$\sumjn a_jx_jb_j\in K$$
for all $x_j\in K$ and $a_j\in A$, $b_j\in B$ satisfying $\sumjn a_ja_j^*\leq1$,
$\sumjn b_j^*b_j\leq1$. 

\begin{theorem}\label{obip}\cite[3.8, 3.9]{M7} Let $K$ be an $A,B$-absolutely convex
subset of a bimodule $X\in\asob$.  If $K$ is closed in the $A,B$-topology, then 
for each $x\in X\setminus K$ there exist normal cyclic Hilbert modules
$\h$ over $A$ and $\k$ over $B$ and a map $\phi\in\acbbb{X}{\bkh}$ 
such that $\norm{\phi(y)}\leq 1$ for all $y\in K$
and $\norm{\phi(x)}>1$. If $X\in\andob$ and $K$ is weak* closed then
$\phi$ can be chosen weak* continuous.
\end{theorem}

For bimodules $U\in\aob$ and $V\in{_B{\rm OM}_A}$ we denote by $U{_A\opr_B}V$ the quotient
of the maximal operator space tensor product $U\opr V$ by the closed subspace
$N$ generated by $\{aub\otimes v-u\otimes bva:\ a\in A,\ b\in B,\ u\in U,\ v\in V\}$. Consider 
the natural
completely isometric isomorphism (\cite[(1.51)]{BLM}, \cite[4.1]{Pi})
$$\iota:\cbb{U}{\du{V}}\to \du{(U\opr V)},\ \ \iota(\phi)(u\otimes v)=\phi(u)(v)\ 
\ \ (\phi\in\cbb{U}{\du{V}})$$
and note that $\iota(\phi)$ annihilates $N$ if and only if $\phi\in
\acbbb{U}{\du{V}}$, where $\du{V}$ is the dual $A,B$-bimodule of $V$ 
(Definition \ref{db}). Thus, we have completely isometrically 
\begin{equation}\label{94} \acbbb{U}{\du{V}}=\du{(U{_A\opr_B}V)}.
\end{equation}

Now we turn to the definition of bimodule duality.
\begin{definition}\label{md} Given operator algebras $A\subseteq\bh$ and 
$B\subseteq\bk$ (containing the identity operators), the {\em bimodule dual} 
(with respect to $\h$ and $\k$) of a bimodule $X\in\aob$ is the 
$\com{A},\com{B}$-bimodule
$\md{X}=\acbbb{X}{\bkh}$, where  
$$(\com{a}\phi\com{b})(x):=\com{a}\phi(x)\com{b}\ \ (\phi\in\md{X}).$$
If $\h$ and $\k$ are proper, 
we emphasize this by writing $\pmd{X}$ instead of $\md{X}$ and call $\pmd{X}$ a 
{\em proper bimodule dual} of $X$.
\end{definition}

We could replace in the above definition the inclusions $A\subseteq\bh$ and
$B\subseteq\bk$ by more general (normal) representations without essentially
changing the ideas, only the notation would be more complicated. If $\h$ and
$\k$ are separable and proper then, as we already remarked, they are unique
up to a unitary equivalence of modules and consequently the proper duals are essentially
unique in this case.
From (\ref{94}) (and using Theorem \ref{thbez})
we deduce by standard arguments:

\begin{proposition}\label{pn} $\md{X}$ is a normal dual operator
$\com{A},\com{B}$-bimodule if $X\in\aob$.
\end{proposition}

\begin{definition}\label{mp} For von Neumann algebras $A\subseteq\bh$ and 
$B\subseteq\bk$ and a bimodule $X\in\adob$,  the $\com{A},\com{B}$-bimodule
$\mpd{X}=\ancbbb{X}{\bkh}$ is called the
{\em bimodule predual}  of $X$. If $\h$ and $\k$ are proper 
then $\mpd{X}$ is denoted also by $\pmp{X}$.
\end{definition}

The following theorem was proved
in \cite{M9} in the 
case $B=A$, but the same proof works in general. 

\begin{theorem}\label{mpmd} If $X\in\anob$, then
$\pmp{(\pmd{X})}$ is the smallest strong $A,B$-bimodule containing $X$. In particular, 
$\pmp{(\pmd{X})}=X$ if and only if $X$ is strong.
\end{theorem}

Now we recall the definition of the (extended) Haagerup tensor product of modules.
For two modules $X\in{\rm NOM}_B$ and $Y\in{_B{\rm NOM}}$ the completion of
the algebraic tensor product $X\otimes_BY$ with the norm  
$$h(w)=\inf\{\norm{\sumjn x_jx_j^*}^{1/2}\norm{\sumjn y_j^*y_j}^{1/2}:\ w=\sumjn x_j\otimes_B y_j\}$$
is the {\em Haagerup tensor product $X\hg_BY$} \cite{BLM}. A typical element $w\in
X\hg_BY$ can be represented as $w=\sum_{j=1}^{\infty}x_i\otimes_By_i$, where 
the two series $\sum_{j=1}^{\infty}x_jx_j^*$ and
$\sum_{j=1}^{\infty}y_j^*y_j$ are norm convergent. We write
this as 
\begin{equation}\label{hr}w=x\odot_By,
\end{equation}
where $x\in\rowj{X}$, $y\in\colj{Y}$ and $\bj=\{1,2,\ldots\}$.

The
{\em extended Haagerup tensor product $X\ehg_BY$}  consists
of all `formal expressions' (\ref{hr}), where $x\in\rowj{X}$ and $y\in\colj{Y}$
for some (infinite) index set $\bj$.  To explain the 
term `formal expression', 
we may assume (by the CES theorem, \cite[3.3.1]{BLM}) that $X,Y,B\subseteq\bh$ for a Hilbert
space $\h$ and regard  $w=x\odot_By$ as  completely bounded map 
$\com{b}\mapsto x\com{b}y$ from $\com{B}$  into $\bh$. From
\cite[3.2]{M3} we have that 
\begin{equation}\label{200} x\odot_By=0\Longleftrightarrow\exists\ 
\mbox{a projection}\ P\in\matj{B}\ \mbox{such that}\ xP=0\ \mbox{and}\ Py=y.
\end{equation}
Thus, $X\ehg_BY$ is
defined as the space of all maps in $\cbb{\com{B}}{\bh}$ that can be represented
in the form (\ref{hr}) with $x\in\rowj{X}$ and $y\in\colj{Y}$ for some cardinal
$\bj$. (The two sums
$\sumj x_jx_j^*$ and $\sumj y_j^*y_j$ are now weak* convergent.)
If $X\in\asob$ and $Y\in{_B{\rm SOM}_C}$ then $X\ehg_BY\in{_A{\rm SOM}_C}$ and for
each $w\in X\ehg_BY$ 
\begin{equation}\label{ehn}\cbnorm{w}=\inf\{\norm{x}\norm{y}:\ w=x\odot_By,\ \  
x\in\rowj{X},\ y\in\colj{Y}\}.
\end{equation}
For more see \cite{M3} and, for alternative approaches in the case $B=\bc$,
\cite{BS}, \cite{ER3}. 
We shall need the following basic property of the symbol $\odot_B$:
\begin{equation}\label{odotb}xb\odot_By=x\odot_Bby,\ \ 
(b\in\matj{B},\ x\in\rowj{X},\ y\in\colj{Y}).
\end{equation}

\begin{remark}\label{reh=h}Since for Hilbert space vectors
$(\xi_j)\in\colj{\h}$ the sum $\sumj
\norm{\xi_j}^2$ is convergent, for a Hilbert $A$-module $\h$ and $X\in_A{\rm SOM}$
we have that
$$X\ehg_A\h=X\hg_A\h\ \ \mbox{and}\ \ \h^*\ehg_AX=\h^*
\hg_AX.$$
\end{remark}

\section{Basic duality for normal operator bimodules}

{\em In this section $A$, $B$ and $C$
are von Neumann algebras and the bimodule duality is defined using fixed faithful 
normal Hilbert modules
$\h$, $\k$, $\l$ over $A$, $B$, $C$ (resp.).}

\begin{definition}\label{d31} Given $X\in\aob$ and $Y\in{{_B{\rm OM}_C}}$, let
$\mdbxy$ denote the subspace of the $A,C$-bimodule dual of $X\hg_BY$
consisting of all $\Omega\in\md{(X\hg_BY)}$ such that the map 
$B\ni b\mapsto\Omega(x\otimes_Bby)$ is weak* continuous for all $x\in X$,
$y\in Y$.
\end{definition}

A part of the development in this section is based on the following extension
of a result of Blecher and Smith \cite{BS}.

\begin{theorem}\label{t32}  If $X\in\anob$ and $Y\in{_B{\rm NOM}_C}$ then
$\mdbxy=\md{X}\ehg_{\com{B}}\md{Y}$ completely isometrically as 
$\com{A},\com{C}$-bimodules.
\end{theorem}

\begin{proof} Consider the natural map
$\iota:\md{X}\ehg_{\com{B}}\md{Y}\to\mdbxy$ defined by
$$\iota(\phi\odot_{\com{B}}\psi)(x\otimes_By)=\phi(x)\psi(y),$$
where $x\in X$, $y\in Y$, $\phi\in\rowj{\md{X}}$ and $\psi\in\colj{\md{Y}}$.
Note that 
$$\rowj{\md{X}}={\rm CB}_A(X,\rowj{\bkh})_B\subseteq{\rm CB}(X,{\rm B}(\k^{\bj},\h)),$$
hence
$\phi(x)\in\B(\k^{\bj},\h)$ and similarly $\psi(y)\in\B(\l,\k^{\bj})$.
Using (\ref{200}) it can be verified that $\iota$ is a well defined completely 
contractive homomorphism of $\com{A},\com{C}$-bimodules.
To show that $\iota$ is injective, suppose that $\phi\odot_{\com{B}}\psi$ is in
the kernel of $\iota$. This means that 
\begin{equation}\label{31}\phi(X)\psi(Y)=0.
\end{equation}
Since $\psi(Y)$ is a $B$-module, the subspace $[\psi(Y)\l]$ of $\k^{\bj}$ is
invariant under $B$, hence the projection $\com{p}\in\B(\k^{\bj})$ with the 
range $[\psi(Y)\l]$ is
in $\matj{\com{B}}$.
Clearly $\com{p}\psi=\psi$, while (\ref{31}) implies that $\phi\com{p}=0$.
Hence (using 
(\ref{odotb})) $\phi\odot_{\com{B}}\psi=
\phi\odot_{\com{B}}\com{p}\psi=\phi\com{p}\odot_{\com{B}}\psi=0$.

Now, since we have just shown that $\iota$ is injective, it suffices  to prove $\iota$ is a completely quotient map. Let 
$$\Omega\in\matn{\mdbxy}\subseteq{\rm CB}_A(X\hg_BY,\bb{\l^n}{\h^n})_C$$
be a complete contraction.
Then from the well known CSPS
theorem (\cite[p. 17.8]{P}, \cite[1.5.7]{BLM})  
it can be deduced (in the same way as in \cite[proof of 3.9]{M3})
that there exist a normal Hilbert $B$-module $\g$ and
complete contractions $\phi\in\acbbb{X}{\bb{\g}{\h^n}}$ and
$\psi\in{\rm CB}_B(Y,\bb{\l^n}{\g})_C$
such that
$$\Omega(x\otimes_By)=\phi(x)\psi(y)\ \ \ (x\in X,\ y\in Y).$$ 
Since each normal representation of $B$ is contained in a multiple
of the identity representation, we may assume that $\g=\k^{\bj}$ for some
cardinal $\bj$.  Then
$$\phi\in\acbbb{X}{\bb{\g}{\h^n}}=\mat{n}{\bj}{\acbbb{X}{\bkh}}=
{\rm R}_{\bj}(\coln{\md{X}})$$
and 
$$\psi\in{\rm CB}_B(Y,\bb{\l^n}{\g})_C=\mat{\bj}{n}{{\rm CB}_B(Y,\bb{\l}{\k})_C}=
{\rm C}_{\bj}(\rown{\md{Y}}),$$
hence $\phi\odot_B\psi$ is an element of $\coln{\md{X}}\ehg_B
\rown{\md{Y}}=\matn{\md{X}\ehg_B\md{Y}}$ such that  $\norm{\phi\odot_B\psi}\leq 1$ and
$\iota_n(\phi\odot_B\psi)=\Omega$.
\end{proof}

A special case of Theorem \ref{t32} is the following result of Effros and Exel
\cite{EE}.
\begin{corollary}\label{ee}\cite{EE} $\du{(\k^*\hg_{B}\k)}=\com{B}$ if $\k$
is a normal (faithful) Hilbert $B$-module.
\end{corollary}

\begin{proof}We regard $\k$ as a $B,\bc$-bimodule and $\k^*$ as a $\bc,B$-bimodule.
Since $\md{\k}={\rm CB}_B(\k)=\com{B}$
and $\md{(\k^*)}=\com{B}$, we have $\du{(\k^*\hg_B\k)}=\com{B}\ehg_{\com{B}}\com{B}=\com{B}$.
\end{proof}

As an application of Theorem \ref{t32} we can express the bimodule dual of
$X\in\anob$ in terms of usual operator space duality, but first we state a definition.

\begin{definition}\label{d34} If $X\in\aob$, define the Banach 
$\com{B},\com{A}$-bimodule structure on
$\h^*\hg_AX\hg_B\k$  by (using the conventions from Remark \ref{rh})
$$\com{b}(\xi^*\otimes_Ax\otimes_B\eta)\com{a}=\xi^*\com{a}\otimes_Ax\otimes_B
\com{b}\eta.$$
\end{definition} 

Part (i)  of the 
following corollary is known in some form (Na \cite{N}, Blecher \cite{B2}).

\begin{corollary}\label{c35}  For each $X\in\anob$ the following natural maps
are completely isometric isomorphisms of bimodules (and will be regarded as equalities
later on):

(i) $\kappa:\md{X}\to \du{(\h^*\hg_AX\hg_B\k)},\ \ \kappa(\phi)(\xi^*\otimes_A x
\otimes_B\eta)=\inner{\phi(x)\eta}{\xi}.$
Here the $\com{A},\com{B}$-bimodule structure on $\md{X}$ is  as in 
Definition \ref{md}, while the bimodule structure on
$\du{(\h^*\hg_AX\hg_B\k)}$ is dual (in the sense of Definition \ref{db}) to that on 
$\h^*\hg_AX\hg_B\k$ (Definition \ref{d34}).

(ii) $\iota:\h^*\hg_{\com{A}}\md{X}\hg_{\com{B}}\k\to\adb{X},\ \ 
\iota(\xi^*\otimes_{\com{A}}\phi\otimes_{\com{B}}\eta)=\inner{\phi(x)\eta}{\xi}.$
Here the structure of
$B,A$-bimodule on $\h^*\hg_{\com{A}}\md{X}\hg_{\com{B}}\k$ is as in Definition
\ref{d34} (but with $A$ and $B$ replaced by $\com{A}$ and $\com{B}$, resp.),
while $\adb{X}$ inherits its structure from $\du{X}$ (which is 
dual to that on $X$, Definition \ref{db}(i)). In fact,  for each $\rho\in\adb{X}$
there exist a set $\bj$, unit vectors $\xi\in\h^{\bj}$ and $\eta\in\k^{\bj}$
and a map $\phi\in\matj{\md{X}}$ such that $\rho(x)=\inner{\phi(x)\eta}{\xi}$ ($x\in X$) and
$\cbnorm{\phi}=\cbnorm{\rho}$.

(iii) $\du{(\adb{X})}=\mbd{X}.$
\end{corollary}

\begin{proof} The routine verifications that $\kappa$ and $\iota$ are bimodule 
homomorphisms and that the identifications
below are the same as stated in the Corollary will be omitted. 

(i) That
$\kappa$ is a complete isometry follows from Theorem \ref{t32} and the associativity
of the (extended) Haagerup tensor product. Namely, since $\md{\k}=\com{B}$ (as in
the proof of Corollary \ref{ee}) 
and similarly $\md{(\h^*)}=\com{A}$ and $\h$, $\k$ are normal, we have the following complete
isometries (regarded as equalities):
$$\du{(\h^*\hg_AX\hg_B\k)}=\md{(\h^*\hg_AX)}\ehg_{\com{B}}\md{\k}=
\md{(\h^*)}\ehg_{\com{A}}\md{X}=\md{X}.$$

(ii) Regarding $A$ as a  $\bc,A$-bimodule and $B$ as a  $B,\bc$-bimodule,
we have $\md{A}=\h^*$ and $\md{B}=\k$. Thus by Theorem \ref{t32} 
$$\adb{X}=\adb{(A\hg_AX\hg_BB)}=\md{A}\ehg_{\com{A}}\md{X}\ehg_{\com{B}}\md{B}=
\h^*\ehg_{\com{A}}\md{X}\ehg_{\com{B}}\k.$$
The norm equality $\cbnorm{\phi}=\norm{\rho}$ follows from the proof of Theorem \ref{t32}.

(iii) This is an immediate consequence of (i) and (ii).
\end{proof}

\begin{corollary}\label{c36} For each $X\in\anob$ the natural homomorphism
$X\to\mbd{X}$ is completely isometric.
\end{corollary}

\begin{proof} Note that there is a completely contractive projection from $\du{X}$ onto
$\adb{X}$ (see \cite[4.4]{M5} or the proof of Theorem \ref{t59}), hence by Corollary 
\ref{c35}(iii) $\mbd{X}=\du{(\adb{X})}\subseteq \bdu{X}$. 
\end{proof}

The following result is dual to Theorem \ref{mpmd}.

\begin{theorem}\label{t38} For each $X\in\andob$ the natural map
$\iota:X\to\pmd{(\pmp{X})}$ is a completely isometric weak* homeomorphic
isomorphism of $A,B$-bimodules.
\end{theorem}

\begin{proof} Set $Y=\pmp{X}$. It is not hard to verify that $Y$ is a strong 
$\com{A},\com{B}$-subbimodule in $\md{X}$ (see \cite[p. 156]{ER1} if necessary). 
To prove that the natural $A,B$-bimodule complete
contraction
$$\iota:X\to\pmd{Y},\ \ \iota(x)(\phi)=\phi(x)\ \ (\phi\in Y)$$
is completely isometric,
let $x\in\matn{X}$ with $\norm{x}>1$. By Theorem \ref{obip} applied to the normal dual
$\matn{A},\matn{B}$-bimodule $\matn{X}$ (with $K$ the unit ball of $\matn{X}$) there
exist cyclic normal Hilbert modules $\tilde{\g}$ over $\matn{A}$ and 
$\tilde{\l}$ over $\matn{B}$ and a weak* continuous
contractive bimodule map $\tilde{\phi}:\matn{X}\to
\bb{\tilde{\l}}{\tilde{\g}}$ such that $\norm{\tilde{\phi}(x)}>1$. In fact
$\cbnorm{\tilde{\phi}}\leq1$ by a result of Smith quoted below as Theorem \ref{acb}.
An elementary well known argument about Hilbert modules over $\matn{A}$ shows that  $\tilde{\g}=\g^n$ and 
$\tilde{\l}=\l^n$ for some normal Hilbert modules $\g$ over $A$ and $\l$ over $B$
and (since $\tilde{\phi}$ is a homomorphism of 
$\matn{A},\matn{B}$-bimodules) $\tilde{\phi}=\phi_n$,
where  
$\phi\in\ancbbb{X}{\B(\l,\g)}$ (that is, $\tilde{\phi}([x_{ij}])=[\phi(x_{ij})]$ for all
$[x_{ij}]\in\matn{X}$). Since $\tilde{\g}$ and $\tilde{\l}$ are cyclic over
$\matn{A}$ and $\matn{B}$ (resp.), $\g$  and $\l$ are $n$-cyclic over $A$ and $B$
(resp.),
which means that (up to a unitary equivalence) $\g\subseteq\h^n$ and 
$\l\subseteq\k^n$,
where $\h$ and $\k$ are proper  modules used in the definition of
duality. Then $\phi$ may be regarded as an element of $\ancbbb{X}{\matn{\bkh}}=\matn{Y}$ and $\cbnorm{\phi}
\leq1$. Since $\norm{\phi_n(x)}>1$, it follows that $\norm{\iota(x)}>1$ and
$\iota$ must be completely isometric.

Next note that $\iota$ is weak* continuous on the unit ball, hence a weak* homeomorphism
onto the weak* closed subspace $\iota(X)$ in $\pmd{Y}$ by the Krein - Smulian theorem.
Indeed, if $(x_j)$ is a bounded net in $X$ weak* converging to an $x\in X$, 
then for each $\phi\in Y\ (=\pmp{X})$ the net $(\phi(x_j))$ converges
to $\phi(x)$ in the weak* topology of $\bkh$, hence
$$\inner{\iota(x_j)}{\xi^*\otimes_{\com{A}}\phi\otimes_{\com{B}}\eta}=
\inner{\phi(x_j)\eta}{\xi}\to\inner{\phi(x)\eta}{\xi}$$
for all $\xi\in\h$ and $\eta\in\k$. Since elements of the form 
$\xi^*\otimes_{\com{A}}\phi\otimes_{\com{B}}\eta$ generate the predual 
$\h^*\hg_{\com{A}}Y\hg_{\com{B}}\k$ of $\pmd{Y}$, this
proves that $\iota$ is weak* continuous.

Now we may identify $X$ with $\iota(X)$ in $\pmd{Y}$. If $X\ne\pmd{Y}$,
then by Theorem \ref{obip} there exists a nonzero $\phi\in\pmp{(\pmd{Y})}$ 
such that $\phi(X)=0$. But, since $Y$ is a strong $\com{A},\com{B}$-bimodule,
$\pmp{(\pmd{Y})}=Y$ by Theorem \ref{mpmd}. Thus $\phi\in Y$ and therefore
$\phi(X)=0$ implies  $\phi=0$, since $Y=\pmp{X}$. This contradiction proves that
$X=\pmd{Y}$.
\end{proof}

We remark without proof that the restriction to proper duals in Theorem \ref{t38}
is necessary, without this restriction the
map $\iota$ need not be isometric.

\begin{definition}\label{d42} A bimodule $X\in\anob$ is 
{\em $A,B$-reflexive} if the natural complete isometry $X\to\mbd{X}$
is surjective.
\end{definition}

Here is a generalization of the classical characterization of reflexivity.

\begin{proposition}\label{p43} A bimodule $X\in\anob$ is $A,B$-reflexive if and only if
its unit ball $B_X$ is compact in the topology induced by
$\adb{X}$.
\end{proposition}

\begin{proof} By Corollary \ref{c35}(i) $\mbd{X}=\du{(\adb{X})}$. By classical
arguments the unit ball of $\du{(\adb{X})}$ is compact in the topology induced
by $\adb{X}$, with $B_X$ a dense subset.
\end{proof}

As an immediate application of  Proposition \ref{p43} one can deduce that strong
subbimodules of $A,B$-reflexive normal operator 
bimodules are $A,B$-reflexive and that if $X\in\anob$ is $C,D$-reflexive for some von 
Neumann subalgebras $C\subseteq A$ and $D\subseteq B$, then $X$ is $A,B$-reflexive
(since the topology induced by $\adb{X}$ is weaker than that by $X^{_C\sharp_D}$).

Now we consider the (non) reflexivity of the basic bimodule $\bkh$.
\begin{example}\label{e44} The bimodule $\bkh$ is $A,B$-reflexive if and 
only if at least one of the algebras $A$ or $B$ is
atomic and finite. 

To prove this, note that by Proposition \ref{p43} the $A,B$-reflexivity does not depend
on the choice of $\h$ and $\k$ in the definition of duality, hence we may assume
that $A\subseteq\bh$ and $B\subseteq\bk$ are in the standard form, so of the same
type as $\com{A}$ and $\com{B}$, respectively.

If, say $\com{B}$, is atomic and finite then by \cite[3.4]{M9}
$$\md{\bkh}=\acbb{\bkh}=\ancbb{\bkh}=\com{A}\ehg\com{B}.$$ Since the unit ball of
$\matn{\com{A}\hg\com{B}}$
is dense in that of $\matn{\com{A}\ehg\com{B}}$ in the $\com{A},\com{B}$-topology 
(which can be shown by approximating elements of $\matn{\com{A}\ehg\com{B}}=
\coln{\com{A}}\ehg\rown{\com{B}}$ by finite sums similarly as in \cite[p. 33]{M5}), 
it follows that
\begin{equation}\label{444}
{\rm CB}_{\com{A}}(\com{A}\ehg\com{B},\bkh)_{\com{B}}=
{\rm CB}_{\com{A}}(\com{A}\hg\com{B},\bkh)_{\com{B}}=\bkh,
\end{equation}
hence $\mbd{\bkh}=\bkh$ and $\bkh$ is $A,B$-reflexive.

On the other hand, by \cite{EK}
$\md{\bkh}=\acbb{\bkh}=\com{A}\nhg\com{B}=:V$, and $V$
contains $\ancbb{\bkh}=\com{A}\ehg\com{B}=:U$. Now $\md{V}=\mbd{\bkh}$ and
$\md{U}={\rm CB}_{\com{A}}(\com{A}\hg\com{B},\bkh)_{\com{B}}=\bkh$. If
$\mbd{\bkh}=\bkh$, the two strong $\com{A},\com{B}$-bimodules $U$ and $V$ have the same bimodule
dual $\bkh$, hence $U=V$ by Theorem \ref{mpmd}. It follows that $C\nhg D=
{\rm CB}_{\com{C}}(\bb{\k}{\h})_{\com{D}}={\rm NCB}_{\com{C}}(\bb{\k}{\h})_{\com{D}}=C\ehg D$ for all
von Neumann algebras $C\subseteq\com{A}$ and $D\subseteq\com{B}$. If neither
$\com{A}$ nor $\com{B}$ is atomic and finite, we can choose $C$ and $D$ both isomorphic
to $C=L_{\infty}[0,1]$. But, with this choice, $C\nhg C\ne C\ehg C$ since
there exist non-normal completely bounded $C$-bimodule maps on $\B(L_2[0,1])$.
\end{example}

Now we are going to consider very briefly the adjoints of bimodule maps. Again,
for maps between strong bimodules the results resemble the classical ones, but
there is a difference (Proposition \ref{p311}(iii) below).

The {\em bimodule adjoint} of a map $T\in\acbbb{X}{Y}$
is  the $\com{A},\com{B}$-bimodule map 
\begin{equation}\label{ma}\md{T}:\md{Y}\to\md{X},\ \ \md{T}(\psi)=\psi\circ T\ \ 
(\psi\in\md{Y}). 
\end{equation}
If $\md{X}$ and $\md{Y}$ are proper bimodule duals, we write $\pmd{T}$ instead of $\md{T}$. 

The following proposition can be deduced using Theorem \ref{mpmd} by standard
arguments, so we omit its proof.
\begin{proposition}\label{p310} If  $X,Y\in\anob$ with $Y$ strong and 
$T\in{\rm NCB}_{\com{A}}(\pmd{Y},\pmd{X})_{\com{B}}$, then there exists a unique
$S\in\acbbb{X}{Y}$ such that $T=\pmd{S}$.
\end{proposition}

Note that $\cbnorm{\md{T}}\leq\cbnorm{T}$. If $X,Y\in\anob$ and 
$T\in\acbbb{X}{Y}$ then,
using the identification $\md{X}=(\h^*\hg_AX\hg_B\k)^{\sharp}$ from Corollary \ref{c35}(i), $\md{T}$ can be expressed as the usual adjoint of
another completely bounded map $T_h$:
\begin{equation}\label{ma1}\md{T}=\du{T_h},\ \mbox{where}\ 
T_h=1_{\h^*}\otimes_AT\otimes_B1_{\k}:\h^*\hg_AX\hg_B\k\to\h^*\hg_AY\hg_B\k.
\end{equation}

\begin{proposition}\label{p311} Let $X,Y\in\asob$ and $T\in\acbbb{X}{Y}$. Then:

(i) $\cbnorm{\md{T}}=\cbnorm{T}$;

(ii) $T$ is a complete isometry if and only if $\md{T}$ is a completely quotient
map. 

(iii) $\pmd{T}$ is a complete isometry if and only if for each $n\in\bn$ the image 
$T(B_{\matn{X}})$
of the unit ball of $\matn{X}$ is dense in  $B_{\matn{Y}}$ in
the $A,B$-topology.

(iv) If $\pmd{T}$ is a complete isometry and $T$ is injective, then $T$ is a
completely isometric surjection.
\end{proposition}

\begin{proof} Parts (i), (ii) and (iii) can be deduced by classical reasoning
using the operator bipolar Theorem \ref{obip} and Corollary \ref{c36}. 
(Alternatively, using (\ref{ma1}), (i) and (ii) can also be deduced from the corresponding properties 
of the usual completely bounded adjoint operators, but we omit the details.)
To prove (iv), observe that since $T$ is injective, the same holds for
$T_h$ in (\ref{ma1}). (Indeed, if $T_j:X_j\to Y_j$ are injective bimodule maps 
then,
using (\ref{200}),  $T_1\otimes_BT_2:X_1\ehg_BX_2\to Y_1\ehg_BY_2$ can
easily be proved to be injective.)
Then, by classical duality and (\ref{ma1}) $\pmd{T}$ has dense range. On the other
hand, since
$\pmd{T}$ is a weak* continuous isometry and the ball $B_{\pmd{Y}}$ is weak* 
compact,  $B_{\pmd{T}(\pmd{Y})}=
\pmd{T}(B_{\pmd{Y}})$ must be weak* compact. Now the Krein - Smulian theorem
shows that $\pmd{T}(\pmd{Y})$ is weak* closed, hence it follows that $\pmd{T}$ is
surjective. Thus $\pmd{T}$ is a weak* homeomorphism
of the unit balls, hence $R:=(\pmd{T})^{-1}$ is weak* continuous by the Krein - Smulian
theorem. By Proposition \ref{p310} there exists an $S\in\acbbb{Y}{X}$ such that
$R=\pmd{S}$. From $\pmd{S}=(\pmd{T})^{-1}$ we conclude that $T=S^{-1}$; moreover, 
since $S$ and $T$
are complete contractions, both must be completely isometric.
\end{proof}

We conclude this section with some applications to Hilbert W$^*$-modules and correspondences. 
These will not be needed later in the paper. 
Basic facts
about such modules can be found in many sources (e.g.  \cite{BLM}
or \cite{R}). We only recall that if $B$ and $C$ are von Neumann algebras, a 
W$^*$-correspondence from $B$ to $C$ is a self-dual right Hilbert C$^*$-module $F$
over $C$ together with a normal representation of $B$ in the von Neumann algebra $\Ll(F)$ of
all adjointable operators on $F$, hence $F\in{_B{\rm NDOM}_C}$.
In this case Theorem \ref{t32} can be slightly improved.

\begin{proposition}\label{ph1} If $F$ is a W$^*$-correspondence from $B$ to $C$,
then $$(X\hg_BF)^{\natural_{B{\rm nor}}}=(X\hg_BF)^{\natural}=\md{X}\ehg_{\com{B}}
\md{F}$$ for each $X\in\anob$.
\end{proposition}

\begin{proof} By Theorem \ref{t32} it suffices to prove the first equality. 
We have to show that for each $\theta\in(X\hg_BF)^{\natural}$ and $x\in X$
the map $B\ni b\mapsto\theta(x\otimes_Bby)\in\bb{\l}{\h}$ is normal. Consider the
$C$-module map $\theta_x:F\to\bb{\l}{\h}$, $\theta_x(y)=\theta(x\otimes_By)$.
We recall that $F$ is an orthogonally  complemented submodule in
$\colj{C}$ for some cardinal $\bj$ \cite[8.5.25]{BLM}, hence $\Ll(F)$ can be regarded as a
w*-closed self-adjoint subalgebra in $\Ll(\colj{C})=\matj{C}$. 
Extending $\theta_x$ to a map 
$\sigma\in\cbb{\colj{C}}{\bb{\l}{\h}}_C=\rowj{\bb{\l}{\h}}$
(the last equality can be proved by using (\ref{94}), or see the proof
in \cite[5.1]{M9}), $\sigma$ is of the form $\sigma(y)=Ty$ for some
$T\in\rowj{\bb{\l}{\h}}$. It follows that $\theta(x,by)=\sigma(by)=
Tby$, which is w*-continuous in $b$.
\end{proof}

\begin{corollary}\label{ch1} If $X\in\asob$ and $F$ is a W$^*$-correspondence
from $B$ to $C$, then $(X\ehg_BF)^{\natural}=\md{X}\ehg_{\com{B}}\md{F}$.
\end{corollary}

\begin{proof}[Sketch of the proof] First note that the unit ball of 
$\matn{X\hg_BF}$ is dense in the $\bc,C$-topology in unit ball of $\matn{X\ehg_BF}=
\coln{X}\ehg_B\rown{F}$ for each $n$.
(This follows by the argument from \cite[p. 33]{M5}, using the polar 
decomposition of elements $y$ in a Hilbert C$^*$-module of the form $\colj{\rown{F}}$,
with $|y|\in\matn{C}$.) By automatic continuity of $C$-module maps it follows that
$(X\ehg_BF)^{\natural}=(X\hg_BF)^{\natural}$ and now Proposition \ref{ph1} concludes
the proof.
\end{proof}

We shall study the bimodule dual of $X\ehg_BY$ in greater generality in the next section.
Here we note that by \cite[3.1]{B5} for W$^*$-correspondences,  $E\ehg_BF$ is equal to the usual
(self-dual) tensor product $E\overline{\otimes}_BF$. Thus, Corollary \ref{ch1} implies
that tensor product of correspondences behaves nicely under the bimodule duality,
which is observed also in \cite{MS}. However, in \cite{MS} the duality 
is defined in a different way, but we shall show in the following example that the
two ways are equivalent.

\begin{example}\label{eh1} To compute the bimodule dual of a W$^*$-correspondence $E$ from
$A$ to $B$ we use Corollary \ref{c35}(i), (\ref{94}) and the well known equality
$\h^*\hg X=\h^*\opr X$ \cite[1.5.14]{BLM} to get
$$\md{E}=\du{(\h^*\hg_AE\hg_B\k)}=\bb{\h^*}{(E\hg_B\k)^*)}_A.$$
The latter space can be naturally identified with ${\rm B}_A(\h,E\hg_B\k)$ (see (\ref{hil})), which is
essentially the definition of the dual in \cite[3.1]{MS}. Note that ${\rm B}_A(\h,E\hg_B\k)$
is a W$^*$-correspondence from $\com{B}$ to $\com{A}$ for the $\com{A}$-valued inner product
$\inner{x}{y}_{\com{A}}=x^*y$ and the $\com{B}$-module action $\com{b}(x\otimes_B\eta)=
x\otimes_B\com{b}\eta$ \cite{R}. In the special case when $A=\bc=\h$, we have that
$\md{E}\cong E\hg_B\k$, hence 
$$\mbd{E}\cong{\rm CB}_{\com{B}}(E\opr_B\k,\k)=(\k^*\opr_{\com{B}}(E\opr_B\k))^{\sharp}
=(E\opr_BB_{\sharp})^{\sharp}=\cbb{E}{B}_B=E$$
since $E$ is self-dual.
By the comment following Proposition \ref{p43} this shows that every W$^*$-correspondence
is reflexive as an operator bimodule. 
\end{example}

\section{The bimodule dual of the extended Haagerup tensor product of bimodules}

Due to the important role of the extended module Haagerup tensor product, it
is worthwhile to compute its bimodule dual.
Effros and Ruan \cite{ER3} defined the {\em normal Haagerup tensor product} of 
dual operator spaces by $\du{U}\nhg\du{V}:=\du{(U\ehg V)}$.
Using that each bimodule $X\in\andob$ is of the form $X=\pmd{(\pmp{X})}$ by 
Theorem \ref{t38}, we may define the module version of this product.

\begin{definition}\label{d91} For $X\in\andob$ and $Y\in{_B{\rm NDOM}_C}$ let
$$X\nhg_BY=\pmd{(\pmp{X}\ehg_{\com{B}}\pmp{Y})},$$
where $\pmp{X}\ehg_{\com{B}}\pmp{Y}$ is regarded as an 
$\com{A},\com{C}$-bimodule.
\end{definition}

The bimodule $X\nhg_BY$ can be described in the following way, which shows in particular that,
as an operator space, $X\nhg_BY$ is independent of $A$ and $C$.
\begin{theorem}\label{t92} $X\nhg_BY=(X\nhg Y)/N$, where $N$ is the weak* 
closed subspace of $X\nhg Y$ generated by all elements of the form
$xb\otimes y-x\otimes by$ ($x\in X,\ y\in Y,\ b\in B$).
\end{theorem}

\begin{proof} Let $\h$, $\k$ and $\l$ be proper Hilbert modules over $A$, $B$ and
$C$ (resp.) in terms of which the duals are defined. By Corollary \ref{c35}(i)
\begin{equation}\label{201}X\nhg_BY=\pmd{(\pmp{X}\ehg_{\com{B}}\pmp{Y})}=\du{(\h^*\ehg_{\com{A}}
\pmp{X}\ehg_{\com{B}}\pmp{Y}\ehg_{\com{C}}\l)}.
\end{equation}
By \cite{BS} and \cite{ER2}
$\k\ehg\k^*=
\du{(\k^*\hg\k)}=\bk$, hence $\com{B}\subseteq\k\ehg\k^*$ and  
\begin{equation}\label{202}U:=\h^*\ehg_{\com{A}}\pmp{X}\ehg_{\com{B}}\pmp{Y}\ehg_{\com{C}}\l=
\h^*\ehg_{\com{A}}\pmp{X}\ehg_{\com{B}}\com{B}\ehg_{\com{B}}\pmp{Y}
\ehg_{\com{C}}\l
\end{equation}
is an operator subspace of
$$\begin{array}{lll}V:&=&\h^*\ehg_{\com{A}}\pmp{X}\ehg_{\com{B}}\bk\ehg_{\com{B}}
\pmp{Y}\ehg_{\com{C}}\l\\
&=&\h^*\ehg_{\com{A}}\pmp{X}\ehg_{\com{B}}\k\ehg\k^*\ehg_{\com{B}}\pmp{Y}
\ehg_{\com{C}}\l.
\end{array}$$
Note that $X=\du{(\h^*\hg_{\com{A}}\pmp{X}\hg_{\com{B}}\k)}$ 
for each $X\in\andob$. (Namely, by Theorem \ref{t38} $X=\pmd{(\pmp{X})}$; 
now apply Corollary \ref{c35}(i).) It follows that 
$V=\pd{X}\ehg\pd{Y}$. From (\ref{201}) and (\ref{202}) we have  $X\nhg_BY=\du{U}$.
The adjoint of the inclusion $U\to V$ is the weak* continuous completely quotient
map 
$$q:X\nhg Y=\du{V}\to\du{U}=X\nhg_BY$$
with $\ker{q}=\ort{U}$, the annihilator of $U$ in $\du{V}$. It remains to prove
that $\ort{U}=N$ or equivalently, since $N$ is weak* closed, that
$U=N_{\perp}\ (\subseteq V).$

A general element $v$ of $V$ has the form $v=\xi^*\odot_{\com{A}}\phi
\odot_{\com{B}}T\odot_{\com{B}}\psi\odot_{\com{C}}\eta$, where
$$\xi\in\colj{\h},\ \eta\in\colj{\l},\ \phi=[\phi_{ij}]\in\matj{\pmp{X}},\ 
\psi=[\psi_{ij}]\in\matj{\pmp{Y}},\  T\in\matj{\bk}$$
for some cardinal $\bj$. We have that $v\in N_{\perp}$ if and only if
$$\inner{v}{xb\otimes y-x\otimes by}=0\ \ \mbox{for all}\ x\in X,\ y\in Y,\ b\in B.$$
This can be written as
$\inner{(\phi(xb)T\psi(y)-\phi(x)T\psi(by))\eta}{\xi}=0$ or 
\begin{equation}\label{91} \inner{\phi(X)(bT-Tb)\psi(Y)\eta}{\xi}=0.
\end{equation}
Since $[\psi(Y)\eta]$ is a $B$-submodule of $\colj{\k}=\k^{\bj}$, we have
that  $[\psi(Y)\eta]=\com{q}\k^{\bj}$ for
a projection $\com{q}\in\matj{\com{B}}$. Similarly $[\phi(X)^*\xi]=
\com{p}\k^{\bj}$ for some
projection $\com{p}\in\matj{\com{B}}$, and (\ref{91}) is equivalent to the
requirement that
$\com{p}(bT-Tb)\com{q}=0$ for all $b\in B$ or
\begin{equation}\label{92} \com{p}T\com{q}\in\matj{\com{B}}.
\end{equation}

Let $\com{e}\in\matj{\com{A}}$ and $\com{f}\in\matj{\com{C}}$ be the
projections with ranges $[A\xi]$ and $[C\eta]$ (resp.). From $\com{q}\psi(y)
\eta=\psi(y)\eta$ ($y\in Y$) we have that $q^{\prime\perp}[\psi(Y)C\eta]=
q^{\prime\perp}[\psi(Y)\eta]=0$
(since $\psi$ is a $C$-module map), hence $q^{\prime\perp}\psi(Y)\com{f}=0$.
This means that 
\begin{equation}\label{93} q^{\prime\perp}\psi\com{f}=0;\ \ \ 
\mbox{similarly}\ \ \com{e}\phi p^{\prime\perp}=0.
\end{equation}
Finally, it follows that
$$\begin{array}{lll}
v&=&\xi^*\odot_{\com{A}}\phi\odot_{\com{B}}T\odot_{\com{B}}\psi\odot_{\com{C}}\eta\\
 &=&(\com{e}\xi)^*\odot_{\com{A}}\phi\odot_{\com{B}}T\odot_{\com{B}}\psi\odot_{\com{B}}
 \com{f}\eta\\
 &=&\xi^*\odot_{\com{A}}\com{e}\phi\odot_{\com{B}}T\odot_{\com{B}}\psi\com{f}
 \odot_{\com{C}}\eta\\
 &=&\xi^*\odot_{\com{A}}\com{e}\phi\com{p}\odot_{\com{B}}T\odot_{\com{B}}
 \com{q}\psi\com{f}\odot_{\com{C}}\eta\ \ \mbox{(by (\ref{93}))}\\
 &=&\xi^*\odot_{\com{A}}\com{e}\phi\odot_{\com{B}}\com{p}T\com{q}\odot_{\com{B}}
 \psi\com{f}\odot_{\com{C}}\eta\\
 &\in&U\ \mbox{(by (\ref{92}) and (\ref{202}))}.
 \end{array}$$
This (reversible) computation proves that $U=N_{\perp}$.
\end{proof}

If $A$ and $C$ are von Neumann algebras on a Hilbert space $\h$, the space $A\nhg C$ was 
identified by Effros and Kishimoto in \cite{EK} with ${\rm CB}_{\com{A}}(\bh)_{\com{C}}$. If $B$ is
a common von Neumann subalgebra in $A$ and $C$, we have a similar identification
for $A\nhg_BC$ (Theorem \ref{t94}), but for this we first need to extend
a result from \cite[p. 131]{BS}. 

\begin{proposition}\label{p93} Given von Neumann algebras  $T\subseteq{\rm B}(\h_T)$, $A$, $B$,
Hilbert spaces $\h$, $\k$ and normal representations $A\stackrel{\alpha}{\to}\bh$, 
$A\stackrel{\pi}{\to} T$,
$B\stackrel{\beta}{\to}\bk$, $B\stackrel{\sigma}{\to}T$, we have
\begin{equation}\label{203} \B_A(\h_T,\h)\ehg_{\com{T}}\B_B(\k,\h_T)=\ancbbb{T}{\B(\k,\h)}
\end{equation}
completely isometrically by letting each $\com{a}\odot_{\com{T}}\com{b}$ to
act on $T$ as
$$(\com{a}\odot_{\com{T}}\com{b})(t)=\com{a}t\com{b}\ \ (t\in T,\ \com{a}\in\rowj{
\B_A(\h_T,\h)},\ \com{b}\in\colj{\B_B(\k,\h_T)}.$$
\end{proposition}

\begin{proof} It is perhaps well known (and easy) that for two Hilbert spaces 
$\g=\coli{\bc}$ and $\l^*=\rowj{\bc}$ and any operator space $X$ we have
$\coli{\bc}\ehg X\ehg\rowj{\bc}=\mat{\bi}{\bj}{X}=\mat{\bi}{\bj}{\ncbb{\du{X}}{\bc}}
=\ncbb{\du{X}}{\mat{\bi}{\bj}{\bc}}$, hence
\begin{equation}\label{400} \g\ehg X\ehg\l^*=\ncbb{\du{X}}{{\rm B}(\l,\g)}
\end{equation}
completely isometrically. 

In the  case $A=\bc=B$ the proof of the Proposition consists of the following
computation:
$$\begin{array}{lll}
\B(\h_T,\h)\ehg_{\com{T}}\B(\k,\h_T)&=&(\h\ehg\h_T^*)\ehg_{\com{T}}
(\h_T\ehg\k^*)\\
 &=&\h\ehg(\h_T^*\ehg_{\com{T}}\h_T)\ehg\k^*\\
 &=&\h\ehg\pd{T}\ehg\k^*\ \ \mbox{(by Corollary \ref{ee})}\\
 &=&\ncbb{T}{\B(\k,\h)}\ \mbox{(by (\ref{400}))}.
 \end{array}$$

In general, we
have now only to show that each $\theta\in\ancbbb{T}{\B(\k,\h)}$, 
just proved to be of the form $\theta=\com{a}\odot_{\com{T}}\com{b}$ for some 
$\com{a}\in\rowj{\B(\h_T,\h)}$ and $\com{b}\in\colj{\B(\k,\h_T)}$, has 
this form with the addition that
$\com{a}\in\rowj{\B_A(\h_T,\h)}$ and $\com{b}\in\colj{\B_B(\k,\h_T)}$; for
this see the proof of \cite[1.2]{M1}.
\end{proof}
In (\ref{203}) $E:=\B_A(\h_T,\h)$ and $F:=\B_B(\k,\h_T)$ are right
Hilbert W$^*$-modules over $\com{\pi(A)}$ and $\com{\beta(B)}$ (resp.), but the tensor
product is over $\com{T}$ (not over $\com{\pi(A)}$). In the special case, when $A=T$, $\pi={\rm id}$ and $\beta$ 
is the inclusion, (\ref{203}) can be interpreted as $E\overline{\otimes}_{\com{T}}F=
{\rm CB}_{\com{T}}(E^*,F)$, a result of Denizeau and Havet as
stated in \cite[3.3]{B5}. Since this will not be needed here, we shall not explain it further.

The following theorem is a generalization of \cite[2.5]{EK}. 
\begin{theorem}\label{t94} Let $B\to A\subseteq\bh$ and $B\to C\subseteq\bl$ be 
normal $*$-homo\-morph\-isms of von Neumann algebras 
(so that $A$ and $C$ are $B$-bimodules). Then
$$A\nhg_BC={\rm CB}_{\com{A}}(\B_B(\l,\h),\B(\l,\h))_{\com{C}}.$$
More precisely, the equality here means the completely isometric weak* homeomorphism of $A,C$-bimodules that sends
$a\otimes_Bc$ to the map $x\mapsto axc$ ($x\in\B_B(\l,\h)$).
\end{theorem}

\begin{proof} Let $\k$ be a proper Hilbert $B$-module. 
Regarding $A$ as a $\bc,B$-bimodule and $C$ as a $B,\bc$-bimodule, we have
as special cases of Proposition \ref{p93}: 
$$\pmp{A}={\rm NCB}(A,\k^*)_B=\h^*\ehg_{\com{A}}\B_B(\k,\h),\  
\pmp{C}={\rm NCB}_B(C,\k)=\B_B(\l,\k)\ehg_{\com{C}}\l,$$
hence
$$A\nhg_BC=\du{(\pmp{A}\ehg_{\com{B}}\pmp{C})}=\left(\h^*\ehg_{\com{A}}\B_B(\k,\h)
\ehg_{\com{B}}\B_B(\l,\k)\ehg_{\com{C}}\l\right)^{\sharp}.$$
Since by Proposition \ref{p93} $$\B_B(\k,\h)\ehg_{\com{B}}\B_B(\l,\k)={\rm NCB}_B
(B,\B(\l,\h))_B=\B_B(\l,\h),$$ it follows (by using Remark \ref{reh=h}, the
commutativity and associativity of $\opr$, the identities $\h^*\hg V=\h^*\opr V$,  
$V\hg\l=V\opr\l$ and (\ref{94})) that
$$\begin{array}{lll} A\nhg_BC&=&
(\h^*\ehg_{\com{A}}\B_B(\l,\h)\ehg_{\com{C}}\l)^{\sharp}\\
 &=&(\h^*\opr_{\com{A}}\,\B_B(\l,\h)\opr_{\com{C}}\,\l)^{\sharp}\\
 &\cong&(\B_B(\l,\h)\,{_{\com{A}}\opr_{\com{C}}}\,(\l\opr\h^*))^{\sharp}\\
 &=&{\rm CB}_{\com{A}}(\B_B(\l,\h),\B(\l,\h))_{\com{C}}.
\end{array}$$
\end{proof}

\section{The normal part of an operator bimodule} 

{\em In this section $A\subseteq\bh$ and $B\subseteq\bk$ will be
C$^*$-algebras, $\Phi:A\to{\rm B}(\tilde{\h})$, $\Psi:B\to{\rm B}(\tilde{\k})$ the universal representations and
$\ta=\weakc{\Phi(A)}$ and $\tb=\weakc{\Psi(B)}$ the von Neumann envelopes of
$A$ and $B$, respectively.}

We first recall some basic facts about the universal representation $\Phi$
of a C$^*$-algebra $A$
(see \cite[Section 10.1]{KR} for more details if necessary). Since $\Phi$ is the direct sum of all cyclic representations of $A$ obtained
from the GNS construction,  each 
$\rho\in\du{A}$ is of the form $\rho(a)=\inner{\Phi(a)\eta}{\xi}$ for some vectors
$\xi,\eta\in\tilde{\h}$, therefore $\rho\Phi^{-1}$ has a unique normal extension 
to $\ta$. It follows that $\ta=\bdu{A}$ and that for each $T\in\B(A,\bl)$ the map 
$T\Phi^{-1}$   has a unique weak* continuous
extension $\tilde{T}:\ta\to\bl$. In particular, with $T=i_A:A\to\bh$ the inclusion, 
the map $\tilde{i}_A:\ta\to\overline{A}$ is a normal $*$-homomorphism, hence
\begin{equation}\label{ce} \ker\tilde{i}_A=\ort{P}\ta\ \ \ (\mbox{and similarly}\ \  
\ker\tilde{i}_B=\ort{Q}\tb)
\end{equation}
for some
central projections $P\in\ta$ (and $Q\in\tb$). Finally, recall 
that a map $T\in\B(A,\bl)$ is weak* continuous
if and only if $T(a)=\tilde{T}(P\Phi(a))$ for all $a\in A$.

Now we are going to explain how a dual Banach  $A,B$-bimodule is 
in a canonical way  an $\ta,\tb$-bimodule.

\begin{definition}\label{d54} Given $X=\du{V}\in\adbb$ (as
in Definition \ref{db}), for each $x\in X$ and $v\in V$ let $\oxv\in\du{A}$
and $\rxv\in\du{B}$ be defined by
$$\oxv(a)=\inner{ax}{v}\ \ \mbox{and}\ \ \rxv(b)=\inner{xb}{v}$$
and let $\toxv$ and $\trxv$ be the weak* continuous extensions of $\oxv\Phi^{-1}$ 
and
$\rxv\Psi^{-1}$ to $\ta$ and $\tb$, respectively. Then for  $a\in\ta$, $b\in\tb$ and
$x\in X$ define $ax$ and $xb$ by
\begin{equation}\label{541}\inner{ax}{v}=\toxv(a)\ \ \mbox{and}\ \ \inner{xb}{v}=
\trxv(b).
\end{equation}
It will turn out that this defines an $\ta,\tb$ bimodule structure on $X$, which
will be called {\em the canonical $\ta,\tb$-bimodule}
structure on $X$.
\end{definition}

Relations (\ref{541}) mean that if $a\in\ta$, $b\in\tb$  and $(a_i)$, $(b_j)$ are  
nets in $A$ and $B$ (resp.) 
such that $(\Phi(a_i))$ and $(\Psi(b_j))$ weak* converge to $a$ and $b$  
(resp.), then 
\begin{equation}\label{542} ax=\lim_ia_ix\ \ \mbox{and}\ \ xb=\lim_jxb_j
\end{equation}
in the weak* topology of $X$. 

\begin{remark}\label{r540} Recall (Theorems \ref{om}, \ref{thbez}) that on a dual operator space $X$
each operator left $A$-module structure is given by a $*$-homomorphism $\pi$ 
from $A$ into the von Neumann algebra $A_l(X)$. The above structure of a left $\ta$-module
then necessary comes from the normal extension $\tilde{\pi}:
\ta\to A_l(X)$ of $\pi$. A similar conclusion holds for right modules and $X$ is automatically
a normal dual operator $\ta,\tb$-bimodule.  If $X$ is a general dual Banach bimodule, however, we need to prove that
\begin{equation}\label{asoc} (ax)b=a(xb)\ \ (a\in\ta,\ b\in\tb,\ x\in X).
\end{equation} 
\end{remark}

\begin{proposition}\label{p55} (i) If $X\in\adbb$ then the relations 
(\ref{541}) introduce to $X$ the structure of a  Banach  
$\ta,\tb$-bimodule. Moreover, if $X\in\adob$ then $X$ is a normal dual operator 
$\ta,\tb$-bimodule.

(ii) Each weak* continuous $A,B$-bimodule map $T$ between dual Banach 
$A,B$-bimodules is
automatically an $\ta,\tb$-bimodule map.
\end{proposition}

\begin{proof} (i) The relations $(a_1a_2)x=a_1(a_2x)$ and $x(b_1b_2)=(xb_1)b_2$
($a_k\in\ta,\ b_k\in\tb$) follow easily from (\ref{542}). 
To prove (\ref{asoc}), chose
nets $(a_i)\subseteq A$ and $(b_j)\subseteq B$ so that $(\Phi(a_i))$ and 
$(\Psi(b_j))$ weak* converge to $a\in\ta$ and $b\in\tb$ (resp.). Then, since
the right multiplication by $b_j$ on $X$ is weak* continuous,
$$(ax)b_j=(\lim_ia_ix)b_j=\lim_i(a_ixb_j)=\lim_i(a_i(xb_j))=a(xb_j).$$
Therefore $(ax)b=\lim_j((ax)b_j)=\lim_j(a(xb_j))$ and we would like to show that
this is equal to $a(xb)$ or, equivalently, that
$$\lim_j\inner{a(xb_j)}{v}=\inner{a(xb)}{v}=\tilde{\omega}_{xb,v}(a)$$
for each $v\in V=\pd{X}$. For this, it suffices to show that  
(for $a\in\ta$) the functional $\tb\ni b\mapsto\tilde{\omega}_{xb,v}(a)$
is normal, which in turn is a consequence of weak compactness of
bounded operators from C$^*$-algebras to preduals of von Neumann algebras  \cite{A}. 
Namely, the weak compactness of the operator 
$T:A\to\du{B}$, $T(a)(b)=\theta(a,b)$, where $\theta(a,b)=\omega_{xb,v}(a)=
\rho_{ax,v}(b)$, implies that the left and the right canonical extensions
of $\theta$ to $\ta\times\tb$ agree (see \cite[p. 12]{DF}). This means
that $\tilde{\omega}_{xb,v}(a)=\tilde{\rho}_{ax,v}(b)$, which is a normal functional in
the variable $b\in\tb$. 

If $X\in\adob$ then, as we have noted in Remark \ref{r540}, $X$  is a normal 
dual operator $\ta,\tb$-bimodule.

(ii) This is a consequence of (\ref{542}) and the weak* continuity of $T$.
\end{proof}

\begin{remark}\label{r56} Given $X\in\abb$, $\du{X}$ is a dual Banach
$B,A$-bimodule (in the sense of Definition \ref{db}), hence by Proposition
\ref{p55} $\du{X}$ is 	
canonically a $\tb,\ta$-bimodule. Now on $\bdu{X}$ we have two
$\ta,\tb$-bimodule structures:

(i) The dual $\ta,\tb$-bimodule in the sense of Definition
\ref{db}, that is $\inner{aFb}{\theta}=
\inner{F}{b\theta a}$ ($a\in\ta$, $b\in\tb$, $\theta\in\du{X}$, $F\in\bdu{X}$);
we denote this bimodule by $\bdu{X}_d$.

(ii) The canonical $\ta,\tb$-bimodule as in Definition \ref{d54},
that is $aF=\lim_ia_iF$ and $Fb=\lim_jFb_j$ in the weak* topology of $\bdu{X}$, 
where $(a_i)\subseteq A$ and $(b_j)\subseteq B$ are nets such that  
$(\Phi(a_i))\to a$ and $(\Psi(b_j))\to b$  and where  $\bdu{X}$
(as an $A,B$-bimodule)  is dual to the $B,A$-bimodule
$\du{X}$. 

If $\bdu{X}_d$ is a normal $\ta,\tb$-bimodule, then $\bdu{X}_d=\bdu{X}$ by
continuity since
$\bdu{X}_d$ and $\bdu{X}$ agree as $A,B$-bimodules. 
\end{remark}

\begin{proposition}\label{p57} If $X\in\aob$, then 
$\bdu{X}_d$ is a normal dual operator $\ta,\tb$-bimodule, hence $\bdu{X}_d=\bdu{X}$.
\end{proposition}

\begin{proof} There exist a Hilbert space $\l$, 
representations $\pi:A\to\bl$ and $\sigma:B\to\bl$ and a completely isometric
$A,B$-bimodule embedding $X\subseteq\bl$. Then $\bdu{X}_d\subseteq\bdu{\bl}_d=
\widetilde{\bl}$ (the universal von Neumann envelope of $\bl$), hence it suffices to prove that $\widetilde{\bl}$ is a normal
$\ta,\tb$-bimodule, where 
$$\inner{ax}{\theta}=\inner{x}{\theta a}\ \ \mbox{and}\ \ 
\inner{xb}{\theta}=\inner{x}{b\theta}\ \ (a\in\ta,\ b\in\tb,\ \theta\in\du{\bl},
\ x\in\widetilde{\bl}).$$
Here $\inner{x}{\theta a}$ means $\widetilde{\theta a}(x)$, where 
$\widetilde{\theta a}$ is the normal extension of the functional $\theta a\in\du{\bl}$
to $\widetilde{\bl}$. But, since the multiplication $\ta\times\widetilde{\bl}\ni(a,x)
\mapsto ax$ is separately weak* continuous in both variables (for 
$ax$ is just the internal product $\bdu{\pi}(a)x$ in
$\widetilde{\bl}$ and $\bdu{\pi}:\ta=\bdu{A}\to\bdu{\bl}
=\widetilde{\bl}$ is normal), we have that
$\widetilde{\theta a}(x)=\tilde{\theta}(ax)$, where $\tilde{\theta}$ is the weak*
continuous extension of $\theta\in\du{\bl}$ to  $\widetilde{\bl}$. It
follows that $\inner{x}{\theta a}=\tilde{\theta}(ax)$ and, since the map 
$\ta\ni a\mapsto\tilde{\theta}(ax)$ is weak* continuous, $\widetilde{\bl}$ is a 
normal left $\ta$-module. 
Similarly $\widetilde{\bl}$ is a normal right $\tb$-module. 
The identity $\bdu{X}_d=\bdu{X}$ follows now from Remark \ref{r56}.
\end{proof}

\begin{definition}\label{d51} The {\em normal part} of a bimodule $X\in\aob$, 
denoted by $\nor{X}$, is the norm
closure of $\iota(X)$, where  $\iota:X\to\mbd{X}$ is
the natural complete contraction. 
\end{definition}

The name `normal part' may be justified by the universal property of $\nor{X}$ stated in
part (i) of the following proposition.
\begin{proposition}\label{p53} Let $A$ and $B$ be von Neumann algebras
and $X\in\aob$. Then:

(i) $\nor{X}\in\anob$ and the canonical map $\iota\in
\acbbb{X}{\nor{X}}$ has the following properties: (1) $\cbnorm{\iota}\leq1$
and (2) for each $Y\in\anob$ and $T\in\acbbb{X}{Y}$ there exists a unique
$\nor{T}\in\acbbb{\nor{X}}{Y}$ with $\norm{\nor{T}}\leq\norm{T}$ and
$\nor{T}\iota=T$. Moreover, if $X_0\in\anob$ and a map $\iota_0\in\acbbb{X}{X_0}$
also has the properties (1) and (2) (with $\iota$ replaced by $\iota_0$), then there exists a completely isometric 
$A,B$-bimodule isomorphism
$\sigma:\nor{X}\to X_0$  such that $\iota_0=\sigma\iota$.

(ii) If $Y\in\anob$ and $\phi\in\cbb{X}{Y}$ is weakly $A,B$-continuous in the
sense that $\rho\phi\in\adb{X}$ for each $\rho\in\adb{Y}$,
then there exists a (unique) map $\nor{\phi}\in\cbb{\nor{X}}{Y}$ such that
$\nor{\phi}\iota=\phi$, and we have that  $\cbnorm{\nor{\phi}}=\cbnorm{\phi}$ 
and $\nor{\phi}$ is weakly $A,B$-continuous.
In particular $\adb{X}=\adb{(\nor{X})}$ completely
isometrically.
\end{proposition}

\begin{proof} (i) By Proposition \ref{pn} $\mbd{X}$ (hence also $\nor{X}$) is a 
normal operator $A,B$-bimodule. 
If $\iota_Y:Y\to\mbd{Y}$ is the canonical inclusion (completely
isometric by Corollary \ref{c36} since $Y$ is normal), then $\iota_YT=\mbd{T}\iota_X$, hence we may
simply set $\nor{T}=\mbd{T}|\nor{X}$. The rest of (i) is evident, by
elementary categorical arguments.

(ii) If $\rho\in\du{X}$, then $\rho$ is a weak* continuous functional on the
normal dual operator $\ta,\tb$-bimodule $\bdu{X}$, hence it follows from Corollary
\ref{c35}(ii) 
that there exist an index set $\bj$, unit vectors $\xi\in
\tilde{\h}^{\bj}$, $\eta\in\tilde{\k}^{\bj}$ and a map
$\psi\in{\rm CB}_{\ta}(\bdu{X},{\rm B}(\tilde{\k}^{\bj},\tilde{\h}^{\bj}))_{\tb}$
such that 
$$\rho(x)=\inner{\psi(x)\eta}{\xi}\ \ (x\in\bdu{X})$$
and $\cbnorm{\psi}=\cbnorm{\rho}$. 
If in addition $\rho\in\adb{X}$ then, since the functionals $A\ni a\mapsto\rho(ax)$
and $B\ni b\mapsto\rho(xb)$ are normal, it follows by \cite[10.1.13]{KR}
that 
$$\rho(x)=\rho(PxQ)=\inner{\psi(PxQ)\eta}{\xi}=\inner{P\psi(x)Q\eta}{P\xi}\ \ 
(x\in X).$$
We may regard the map $X\ni x\mapsto P\psi(x)Q$  as an $A,B$-bimodule
map $\psi_0$ from $X$ into the normal operator $A,B$-bimodule 
${\rm B}(Q\tilde{\k}^{\bj},P\tilde{\h}^{\bj})$, hence by part (i) there exists
a map $\psi_1\in\acbbb{\nor{X}}{{\rm B}(Q\tilde{\k}^{\bj},P\tilde{\h}^{\bj}}$
such that $\psi_0=\psi_1\iota$ and $\cbnorm{\psi_1}\leq\cbnorm{\psi}$. 
With $\nor{\rho}\in\adb{(\nor{X})}$ defined by 
$$\nor{\rho}(v)=\inner{\psi_1(v)Q\eta}{P\xi}\ \ \ (v\in\nor{X}),$$
we clearly have that $\rho=\nor{\rho}\iota$ and $\norm{\nor{\rho}}\leq
\cbnorm{\psi_1}\leq\norm{\rho}$. The reverse inequality, $\norm{\rho}\leq
\norm{\nor{\rho}}$, follows from $\rho=\nor{\rho}\iota$ since $\cbnorm{\iota}\leq1$.
Since $\iota(X)$ is dense in $\nor{X}$, $\nor{\rho}$ is unique.

For a more general weakly $A,B$-continuous map $\phi\in\cbb{X}{Y}$, we 
regard $Y$ as a normal operator $A,B$ subbimodule in $\B(\l,\g)$ for some
normal Hilbert modules $\g$ and $\l$ over $A$ and $B$, respectively. Since
$\omega\phi\in\adb{X}$ for each $\omega\in\B(\l,\g)_{\sharp}$, we have from
the previous paragraph that $\omega\phi(\ker\iota)=0$. Thus, $\phi(\ker\iota)
=0$ and therefore there exists a unique map $\nor{\phi}:\nor{X}\to Y$
such that $\phi=\nor{\phi}\iota$. We shall omit the verification that
this map $\nor{\phi}$ satisfies all the requirements. 
\end{proof}
Finally, we can describe the module bidual $\mbd{X}$ and the normal part 
$\nor{X}$ of an operator bimodule $X$ in a useful alternative way.

\begin{theorem}\label{t59} Let $A$, $B$ be von Neumann algebras and $X\in\aob$.
Regard $X$ as an $A,B$-subbimodule in $\bdu{X}$ 
and let $P\in\ta$, $Q\in\tb$ be the central projections as in (\ref{ce}). Then
$\mbd{X}=P\bdu{X}Q$ and $\nor{X}$ is the norm closure of $PXQ$ in $\bdu{X}$.
For $x\in\matm{X}$ (with $\iota:X\to\nor{X}$ the canonical map) we have that 
\begin{equation}\label{59} \norm{\iota_m(x)}=\inf\left(\sup_j\norm{a_jxb_j}\right),
\end{equation} 
where the infimum is taken either over all nets $(a_j)$ and $(b_j)$ in the unit balls of $A$
and $B$ (respectively) that strongly converge to $1$ or over all nets of projections 
$(a_j)\subseteq A$ and $(b_j)\subseteq B$ converging to $1$.
\end{theorem}

\begin{proof} Since $\adb{X}$ consists of all $\rho\in\du{X}$ such that the two
maps $A\ni a\mapsto\rho(ax)$ and $B\ni b\mapsto\rho(xb)$ are normal and
since a functional $\omega$ on $A$ is normal if and only if $\rho=P\rho$ (and
similarly for $B$), it follows that $\adb{X}=Q\du{X}P$. Since the 
$\ta,\tb$-bimodule
$\bdu{X}$ is dual to the $\tb,\ta$-bimodule  $\du{X}$
by Proposition \ref{p57}, this implies that $\du{(\adb{X})}=
P\bdu{X}Q$. By Proposition \ref{p53} we have that $\md{X}=\md{(\nor{X})}$ and
$\adb{X}=\adb{(\nor{X})}$, hence (applying Corollary \ref{c35}(iii) to $\nor{X}$)
\begin{equation}\label{560}\mbd{X}=\mbd{(\nor{X})}=\du{(\adb{(\nor{X})})}
=\du{(\adb{X})}=P\bdu{X}Q.
\end{equation}
Now it
follows from the definition that
$\nor{X}$ is just the norm closure of $PXQ$.

If $(a_j)$ and $(b_j)$ are nets in the unit balls of $A$ and $B$ (resp.) 
converging to 
$1$ in the strong operator topology, then $\norm{\iota_m(x)}=\sup_j\norm{a_j\iota_m(x)b_j}\leq
\sup_j\norm{a_jxb_j}$ since $\nor{X}$ is  normal and $\cbnorm{\iota_m}\leq1$. This
proves the inequality $\leq$ in (\ref{59}). To prove the reverse inequality,
choose nets $(a_j)$ and $(b_j)$ in the unit balls of $A$ and $B$ so that 
$(\Phi(a_j))$ and $(\Psi(b_j))$ strongly converge to
$P$ and $Q$, respectively. (Note that then $(a_j)$ and $(b_j)$ must converge to $1$
since the normal extensions of $\Phi^{-1}$ and $\Psi^{-1}$ are strongly continuous
on bounded sets and
map $P$ and $Q$ to $1$.) Since $\norm{\iota_m{(x)}}=\norm{PxQ}$ and $\bdu{X}$ is a normal
operator $\ta,\tb$-bimodule, 
$$\norm{\iota_m(x)}=\norm{PxQ}=\sup_j\norm{\Phi(a_j)x\Psi(b_j)}=\sup_j\norm{a_jxb_j}.$$
We may replace in this equality each $a_j$ (resp. $b_j$) with the range projection
$R(a_j)\in A$ (resp. $R(b_j)\in B$) since $a_j\leq R(a_j)\leq1$.
\end{proof}

\section{Central  bimodules}

In this section we consider normality for central bimodules. A slightly more general 
version of central bimodules than defined below  is studied also in \cite{B3}.
\begin{definition}\label{d46} A bimodule $X$ over an abelian operator algebra $C$ is
called {\em central} if $cx=xc$ for all $x\in X$ and $c\in C$. The classes of
central $C$-bimodules among, operator and
normal operator bimodules are denoted by  $\coc$ and 
$\cnoc$, respectively.
\end{definition}

\begin{remark}\label{r470} If $C$ is a C$^*$-subalgebra of the center of a
C$^*$-algebra $A$, $J$ is a closed ideal in $C$ and $X\subseteq A$, 
then $d(x,[JA])=d(x,[JX])$ for each $x\in X$, where
$d(x,S)$ denotes the distance of $x$ to a set $S$. This, probably well known
fact, follows by choosing an approximate identity $(e_j)$ for $J$ and noting 
that (since $(e_j)$ is also an approximate identity for $[JA]$) 
$d(x,[JA])=\lim_j\norm{(x-e_jx)}\geq d(x,[JX])$ (see \cite[p. 300]{KR}).
\end{remark}

\begin{remark}\label{r47} For an abelian C$^*$-algebra $C$ we 
denote by $\Delta$ the spectrum of $C$ and
by $C_t$ the kernel of a character $t\in \Delta$. For a bimodule $X\in\coc$ we 
consider the quotients 
$X(t)=X/[C_tX]$. Given $n\in\bn$ and $x\in\matn{X}$ we denote by $x(t)$ the coset 
of $x$ in $\matn{X}/[C_t\matn{X}]$. We shall need to know that
the function
\begin{equation}\label{47}\Delta\ni t\mapsto\norm{x(t)}\end{equation}
is upper semicontinuous and that
\begin{equation}\label{471}\norm{x}=\sup_{t\in\Delta}\norm{x(t)}.\end{equation}
This is known from \cite[p. 37, 41]{DG} and \cite[p. 71]{Po}, but (to avoid
Banach bundles) we provide now a different short argument. We may assume that $X, C\subseteq\bl$ for some Hilbert space $\l$.
Since $X$ is central, $X\subseteq\com{C}$, hence $\matn{X}\subseteq\matn{\com{C}}=:
A$ and $C$ is identified with the center of $A$. Using Remark \ref{r470} we have
that $\norm{x(t)}=d(x,[C_t\matn{X}])=d(x,[C_tA])$, which is just the norm of
the coset of $x$ in $A/[C_tA]$. Now (\ref{471}) and the continuity of the function 
(\ref{47}) follow from  \cite[p. 232]{G}.
We shall call the embedding $$X\to\oplus_{t\in\Delta} X(t),\ \ 
x\mapsto(x(t))_{t\in\Delta}$$
the {\em canonical decomposition} of $X$. 
\end{remark}

{\em Throughout the rest of the section $C$ is an abelian von Neumann
algebra.}

\begin{lemma}\label{l40} A bimodule $X\in\coc$  is normal
if and only if $pX$ is a normal $pC$-bimodule for each $\sigma$-finite projection
$p\in C$. If $C$ is $\sigma$-finite, then $X$ is normal if and only if
\begin{equation}\label{40} \lim_j\norm{p_jx}=\norm{x}
\end{equation}
for each $x\in\matn{X}$ ($n\in\bn$) and each sequence of projections $p_j\in
C$ increasing to $1$. 
\end{lemma}

\begin{proof}  We may assume that $C$ is $\sigma$-finite, for in
general $C$ is a direct sum of $\sigma$-finite subalgebras and $X$ 
(being central)
also decomposes in the corresponding $\ell_{\infty}$-direct sum. Then 
by Theorem \ref{tno} we have to prove   that for each $n\in\bn$, each 
$x\in\matn{X}$ and sequence $(e_j)$ of projections in $\matn{C}$ increasing
to
$1$ the sequence $(\norm{e_jx})$ converges to $\norm{x}$. Suppose the contrary, that
for an $x$ and  a sequence of projections $(e_j)$ we have 
$$\norm{e_jx}\leq M\ \ \mbox{for some constant}\ \ M<\norm{x}.$$
Let $\tau$ be the canonical normal central trace on $\matn{C}$, the values of
which on projections of $\matn{C}$ are of the form $\frac{k}{n}p$, where $p\in C$
is a projection and $k\in\{0,1,\ldots,n\}$. For each $j$ set
$\Delta_j=\{t\in\Delta:\ \tau(e_j)(t)=1\}$, a clopen subset of $\Delta$, and let
$p_j\in C$ be the characteristic
function of $\Delta_j$. Since the sequence $(e_j)$ increases to $1$ and 
$\tau$ is
weak* continuous, $\Lambda:=\bigcup_j\Delta_j$ is dense in $\Delta$. (Otherwise
$\Lambda_0:=\Delta\setminus\overline{\Lambda}$ would be a nonempty open set
such that $\tau(e_j)(t)\leq 1-1/n$ for all $j$, which is impossible since $e_j\nearrow1$.) 
It follows that the 
sequence $(p_j)$
also increases to $1$. For $t\in\Delta_j$, $e_j(t)\in\matn{C}(t)=\matn{\bc}$ is
a projection with the normalized trace equal to $1$, hence $e_j(t)=1$. This implies
that $e_jp_j=p_j$, hence $\norm{p_jx}\leq\norm{e_jx}\leq M<\norm{x}$ 
for all $j$; but this is in contradiction with the assumption (\ref{40}).
\end{proof}

\begin{proposition}\label{p48} A bimodule $X\in\coc$ is normal if and only if for
each $n\in\bn$ and each $x\in\matn{X}$ the function $\Delta\ni t\mapsto \norm{x(t)}$
is continuous.
\end{proposition}

\begin{proof} If $X$ is normal, then  we may 
assume that $X\subseteq\com{C}$, 
the commutant of $C$ in $\bh$  for a  normal Hilbert
$C$-module $\h$, hence $\matn{X}$ is contained in the commutant of $C$ in 
${\rm B}(\h^n)$ and the continuity of (\ref{47}) follows from Remark \ref{r470}
and \cite[p. 233, Lemma 10]{G}.

For the converse, by Lemma \ref{l40} we may assume that $C$ is $\sigma$-finite
and we have to prove the condition (\ref{40}). Let $\Delta_j$ be the clopen
subset of $\Delta$ corresponding to $p_j$, where $p_j$ are projections as in 
Lemma \ref{l40}. Since $p_j\nearrow1$, $\bigcup_j\Delta_j$ is dense in $\Delta$.
Thus, using (\ref{471}), the continuity of the functions $t\mapsto\norm{x(t)}$ implies
that $\norm{x}=\lim_j\sup_{t\in\Delta_j}\norm{x(t)}=\lim_j\norm{p_jx}$.
\end{proof}

\begin{proposition}\label{p49} Let $X\in\cnoc$ be a strong bimodule and 
$Y\subseteq X$ a subbimodule.
Then the quotient $X/Y$ is a normal operator bimodule if and only if $Y$ is
strong and in this case $X/Y$ is also strong. 
\end{proposition}

\begin{proof} It was observed in \cite[p. 204]{M6} that $X/Y$ is normal only if
$Y$ is strong. For the converse, assume that $C$ is $\sigma$-finite
and that the condition in Lemma \ref{l40} for normality is not satisfied.  
Then there exist an  $\dot{x}\in\matn{X/Y}$, a sequence of 
projections $(p_j)$ in $C$ increasing to $1$  and a constant 
$M<\norm{\dot{x}}$ such that $\norm{p_j\dot{x}}< M$ for all $j$. Put $q_0=p_0$
and $q_j=p_j-p_{j-1}$ if $j\geq1$. Let $x\in\matn{X}$ be any representative of
the coset $\dot{x}$. By definition of the quotient norm for each $j$ there
exists an element $y_j=q_jy_j\in\matn{Y}$ such that $\norm{q_jx-y_j}<M$.
Since the sequence $(y_j)$ is bounded and $Y$ is strong, the sum
$$y:=\sum_{j=0}^{\infty}q_jy_j=\sum_{j=0}^{\infty}q_j(y_jq_j)$$
defines an element of $Y$. But then the estimate
$$\norm{x-y}=\norm{\sum_jq_j(x-y)q_j}=\sup_j\norm{q_j(x-y)}\leq\sup_j
\norm{p_j(x-y_j)}<M$$
implies that $\norm{\dot{x}}<M$, which is in contradiction with the choice of
$M$. 

To verify that $X/Y$ is a strong left $C$-module (hence a strong $C$-bimodule
since it is central), let $(p_j)$ be an orthogonal
family of projections in $C$ and $(\dot{x}_j)$ a family of elements in $X/Y$ such
that the sum  $\sum_j\dot{x}_j^*\dot{x}_j$ converges in the strong operator topology
of some $\bh$ containing $X/Y$ as a normal operator $C$-bimodule. 
We can choose for each $\dot{x}_j$ a representative $x_j\in X$ so that the
set $(x_j)_j$ is bounded, and then $x:=\sum_jp_jx_j=\sum_jp_jx_jp_j\in X$.
Since the quotient map $Q:X\to X/Y$ is a bounded $C$-bimodule map
(hence continuous in the $C$-topology), it follows that
$\sum_jp_j\dot{x}_j=\sum_jp_jQ(x_j)=Q(x)$, which shows that $\sum_jp_j\dot{x}_j\in
X/Y$.
\end{proof}

For central bimodules we can now improve Proposition \ref{p311}.

\begin{corollary}\label{c410} If $X,Y\in\cnoc$ are strong and $T\in{\rm CB}_C(X,Y)$,
then $T$ is completely isometric (respectively, completely quotient) if and
only if $\pmd{T}$ is
completely quotient (respectively, completely isometric).
\end{corollary}

\begin{proof} By Proposition \ref{p311} it remains to prove that $T$ is completely
quotient if $\pmd{T}$ is completely isometric. By Proposition \ref{p49}
$X/\ker{T}$ is a strong central $C$-bimodule, hence we consider 
the induced map $\tilde{T}:X/\ker{T}\to Y$. Since $\pmd{\tilde{T}}:\pmd{Y}\to
\pmd{(X/\ker T)}\subseteq\pmd{X}$ is essentially $\pmd{T}$, hence 
completely isometric, and $\tilde{T}$ is injective, it follows
from Proposition \ref{p311}(iv) that $\tilde{T}$ is a completely isometric
surjection, hence $T$ is a completely quotient map.
\end{proof}

\begin{definition}\label{d510} For a function $f:\Delta\to\br$, let $\essup f$ be the infimum of all
$c\in\br$ such that the set $\{t\in\Delta:\ f(t)>c\}$ is meager (= contained
in a countable union of closed sets with empty interiors).

The {\em essential direct sum}, $\eop X(t)$, of a family of Banach spaces 
$(X(t))_{t\in\Delta}$ is defined as the quotient of the $\ell_{\infty}$-direct
sum $\oplus_{t\in\Delta}X(t)$ by the zero space of the seminorm 
$x\mapsto\essup\norm{x(t)}$. Then $\eop X(t)$ with the norm 
$\dot{x}\mapsto\essup\norm{x(t)}$ is a Banach space and we denote by 
$e:\oplus_{t\in\Delta}X(t)\to\eop X(t)$ the quotient map. If $(X(t))_{t\in\Delta}$
is a family of operator spaces, then $\eop X(t)$ is an operator space by the
identification
$$\matn{\eop X(t)}=\eop\matn{X(t)}.$$
\end{definition}

\begin{theorem}\label{t511} Given a bimodule $X\in\coc$ with the canonical decomposition
$\kappa:X\to\oplus_{t\in\Delta}X(t)$, its normal part $\nor{X}$ is just the closure of $e\kappa(X)$ in $\eop X(t)$.
\end{theorem}

\begin{proof} First, to show that $e\kappa(X)$ is a normal operator $C$-module, 
by Lemma \ref{l40} we may assume that $C$ is 
$\sigma$-finite and
it suffices to prove that for each sequence of projections
$p_j\in C$ increasing to $1$ and each $x\in\matn{X}$ the equality
\begin{equation}\label{511} \essup\norm{x(t)}=\lim_j\essup\norm{p_j(t)x(t)}
\end{equation} holds.
With $\Delta_j$ the clopen subset of $\Delta$ corresponding to $p_j$,
$\bigcup_j\Delta_j$ is dense in $\Delta$.  Since the function $\Delta\ni t\mapsto
\norm{x(t)}$ is upper semi-continuous (hence Borel), it agrees outside a meager
set with 
a continuous function $f$ on $\Delta$ by \cite[p. 323]{KR}. Then
$\essup\norm{x(t)}=\sup f(t)$, $\essup\norm{p_j(t)x(t)}=\sup_t p_j(t)f(t)$ and 
$\lim_j\sup_tp_j(t)f(t)=\sup f(t)$ by continuity 
(since $\bigcup_j\Delta_j$ is dense in $\Delta$). 
This implies (\ref{511}).

It remains to show that the closure of $e\kappa(X)$ has the universal property
of $\nor{X}$ from Proposition \ref{p53}(i). Let $Y\in{_C{\rm NOM}_C}$ and 
$T\in{\rm CB}_C(X,Y)$ with $\cbnorm{T}<1$. We have to show that $T$ can be factorized through
$e\kappa(X)$. Replacing $Y$ by the closure of $T(X)$, we may assume that
$Y$ is central. Let $x\in\matn{X}$ and set $y=T_n(x)$. Since $\cbnorm{T}<1$ and $T$
is a $C$-module map, $\norm{y(t)}\leq\norm{x(t)}$ for each $t\in\Delta$.
Set
$$c=\norm{(e\kappa)_n(x)}=\essup_t\norm{x(t)}\ \ \ \mbox{and}\ \ \ V=\{t\in\Delta:\ 
\norm{y(t)}>c\}.$$
Since $Y$ is normal, the function $t\mapsto\norm{y(t)}$ is continuous by Proposition
\ref{p48}, hence $V$ is open. But for each $t\in V$ we have that $c<\norm{y(t)}
\leq\norm{x(t)}$, hence $V$ must be meager by the definition of $c$, hence
$V=\emptyset$ by Baire's
theorem for locally compact spaces. Thus, $\norm{y(t)}\leq c$ for all $t\in\Delta$, 
which means 
that $\norm{T_n(x)}=\norm{y}=\sup_t\norm{y(t)}\leq c=\norm{(e\kappa)_n(x)}$. This
estimate shows that there exists a unique complete contraction
$S:(e\kappa)(X)\to Y$ such that $T=S\circ(e\kappa)$.
\end{proof}

\section{Operator bimodules of a normal representable bimodule}

We begin this section by introducing various classes of Banach bimodules
admitting operator bimodule structures.
\begin{definition}\label{r} (i) \cite{Po} A bimodule $X\in\abb$ is {\em representable}
($X\in\arb$) if for some Hilbert module
$\h$ over $A,B$ there is an isometry in $\B_A(X,\bh)_B$; 
in other words, $X$
can be represented isometrically in $\bh$ as an operator $A,B$-bimodule. 

(ii) If in (i) $A$ and $B$ are von Neumann algebras and $\h$ is
normal over $A$ and $B$, then $X$ is called a {\em normal 
representable} bimodule; the class of all such bimodules is denoted by $\anrb$.

(iii) If $X\in\adbb$ and for some normal Hilbert module
$\h$ over $A$ and $B$ there exists an isometry in $\anbbb{X}{\bh}$, then 
$X$ is called a {\em normal dual 
representable $A,B$-bimodule} ($X\in\andrb$).
\end{definition}

An abstract characterization of normal dual representable bimodules is given
in \cite[4.14]{BM}, but it will not be needed here.

For a representable bimodule $X\in\arb$,  we  define the {\em proper dual} as 
$\pmd{X}=\abbb{X}{\bkh}$, where  $\h$ and $\k$ are  fixed proper modules
over $A$ and $B$, respectively.
Now a bimodule $X\in\aob$ has two proper duals: in the class $\aob$
and in the class $\arb$. But they agree in $\arb$ by the following  result of Smith. 
\begin{theorem}\label{acb} \cite[2.1, 2.2]{S} If $\g$ and $\l$  
are locally cyclic Hilbert modules
over $A$ and $B$ (respectively),  then $\cbnorm{\phi}=\norm{\phi}$ 
for each $\phi\in\B_A(X,\B(\l,\g))_B$ and  $X\in\aob$.
\end{theorem}

By \cite{M6} or \cite{Po} the identities
\begin{equation}\label{22}\m{x}=\sup\norm{axb}\ \ (x\in\matn{X},\ n=1,2,\ldots),
\end{equation}
where the supremum is over all $a$ and $b$ in the unit balls of $\rown{A}$
and $\coln{B}$ (respectively), define on $X$ the {\em minimal operator
$A,B$-bimodule structure}, denoted by $\minab{X}$. 
If $X\in\anrb$, then by considering an isometric representation of $X$ as
a normal subbimodule in some ${\rm B}(\l)$, we see that $\adb{X}$ has enough
functionals to make the natural contraction $\iota:X\to\du{(\adb{X})}$ 
isometric.  Thus, with $Y=\minab{X}$, the completely contractive isometry $\iota
:Y\to\du{(\adb{Y})}$ 
must be completely isometric (otherwise $\iota$ would induce on $X$
an operator $A,B$-bimodule norm structure smaller than $Y=\minab{X}$). Since
$\du{(\adb{Y})}=\pmbd{Y}$ by (\ref{560}) and $\pmbd{Y}$ is normal by Proposition
\ref{pn}, it follows that 
$\minab{X}$ is a normal operator $A,B$-bimodule. Further,
\begin{equation}\label{63} \m{x}=\sup\{\norm{\phi_n(x)}:\ \phi\in\pmd{X},\ 
\norm{\phi}\leq1\}\ \ (x\in\matn{X},\ n=1,2,\ldots).
\end{equation}

\begin{definition}\label{d1} (i) Given $X\in\arb$, the 
{\em maximal operator bimodule norms} are defined
by
\begin{equation}\label{64} \M{x}=\sup\norm{T_n(x)}\ \ (x\in\matn{X},\ n=1,2,\ldots),
\end{equation}
where the supremum is over all  contractions 
$T\in\B_A(X,\matm{\bkh})_B$, with $m\in\bn$ and $\h$, 
$\k$ the Hilbert spaces of the universal representations of $A$ and
$B$, respectively. Denote the operator bimodule so obtained by
$\maxab{X}$.

(ii) If $X\in\anrb$, the {\em maximal normal operator bimodule norms}, denoted
by $\Mn{x}$, are defined by the same formula
(\ref{64}), but with $\h$ and $\k$  (fixed) proper 
Hilbert modules over $A$ and $B$. This operator bimodule is
denoted by $\maxn{X}$.

(iii) If $X\in\andrb$, the {\em maximal normal dual operator bimodule norms}, 
denoted by $\Mnd{x}$, are defined in the same way as $\Mn{x}$, except
that  we now require in addition that the maps $T$ in  (\ref{64}) are weak* continuous.
Denote this operator bimodule by $\maxnd{X}$. 
\end{definition}

Given $X\in\anrb$ 
and $x\in\matn{X}$, since each normal operator $A,B$-bimodule $Y$ is 
contained in a bimodule of the form $\B(\k^{\bj},\h^{\bi})$ with
$\h$ and $\k$ fixed proper modules over $A$ and $B$ (resp.), we 
deduce that $\Mn{x}=\sup\norm{T_n(x)}$, where the supremum is over all contractions
$T\in\B_A(X,Y)_B$ with $Y\in\anob$. We conclude that the operator bimodule $\maxn{X}$
is characterized by the following:
$\maxn{X}$ is a normal operator $A,B$-bimodule and for each $Y\in\anob$
every map $T\in\B_A(X,Y)_B$ is completely bounded from $\maxab{X}$ into
$Y$ with $\cbnorm{T}=\norm{T}$.
There are similar characterizations for $\maxab{X}$ (if $X\in\arb$) and 
$\maxnd{X}$ (if $X\in\andrb$).
From this and the universal property of the
normal part (Proposition \ref{p53}(i)) we deduce:

\begin{corollary}\label{c63} $\maxn{X}$ is the normal part of $\maxab{X}$ 
if $X\in\anrb$.
\end{corollary}

\begin{example} In general ${\rm MAXN}_C(X)_C\ne{\rm MAX}_C(X)_C$ even if $C$ is  
abelian and $X$ is central. To show this, let $U\subseteq V$ be
Banach spaces such that the (completely contractive) inclusion of maximal operator spaces
$\Max{U}\to\Max{V}$ is not completely isometric.
With $\Delta$ the spectrum of $C$ and $t_0\in\Delta$, let
$$X=\{f\in C(\Delta,V):\ f(t_0)\in U\}.$$
We claim that for each $f\in\matn{{\rm MAX}_C(X)_C}$
\begin{equation}\label{65} \norm{f}_{_CM_C}=\max\{\sup_{t\in\Delta}\norm{f(t)}_{\matn{{\rm MAX}(V)}},
\norm{f(t_0)}_{\matn{{\rm MAX}(U)}}\}.
\end{equation}
To show this, it suffices to prove that, when the spaces 
$\matn{X}$ ($n=1,2.\ldots$) are equipped with
the norms defined by the right side of (\ref{65}), each contraction 
$T\in\B_C(X,Y)_C$
into $Y\in{_C{\rm OM}_C}$, is completely contractive. Replacing $Y$ with the closure of
$T(X)$, we may assume that $Y$ is central and therefore has the canonical
decomposition $Y\to\oplus_{t\in\Delta}Y(t)$ (Remark \ref{r47}).
Since $T$ is a $C$-module map, $T$ induces for each $t\in\Delta$ a contraction
$T_t:X(t)\to Y(t)$. Since the operator space
$$X(t)=\left\{\begin{array}{ll}
V&\mbox{if}\ t\ne t_0\\
U&\mbox{if}\ t=t_0
\end{array}
\right.$$
is maximal,  $T_t$ is a complete contraction,
hence so is $T$ (since $\norm{y}=\sup_t\norm{y(t)}$
for each $y\in\matn{Y}$).

Since the inclusion $\Max{U}\to\Max{V}$ is not completely isometric, there exists
$u\in\matn{U}$ with $\norm{u}_{\matn{U}}>\norm{u}_{\matn{V}}$. 
Hence, if $f\in\matn{X}$ is the constant function $f(t)=u$, the function
$t\mapsto\norm{f(t)}$ is not continuous and ${\rm MAX}_C(X)_C$
is not normal by Proposition \ref{p48}. On the other hand, ${\rm MAXN}_C(X)_C$ is 
always normal.
\end{example}

To show that $\maxnd{\cdot}\ne\maxn{\cdot}$, we first need to extend 
\cite[2.8]{B1}.

\begin{proposition}\label{p65} If $X\in\anrb$ then: 
(i) $\pmd{(\maxn{X})}={\rm MIN}_{\com{A}}(\pmd{X})_{\com{B}}$;

(ii) $\pmd{(\minab{X})}={\rm MAXND}_{\com{A}}(\pmd{X})_{\com{B}}$.
\end{proposition}

\begin{proof} (i) Given 
$\phi=[\phi_{ij}]\in\matn{\pmd{(\maxn{X})}}=\B_A(X,\matn{\bkh})_B$,
its norm is $\norm{\phi}=\sup\{\norm{[\phi_{ij}(x)]}:\ x\in X,\ 
\norm{x}\leq1\}$.
Thus, $\pmd{(\maxn{X})}$ is dominated by every  
operator $\com{A},\com{B}$-bimodule norm structure $Z$ on $\pmd{X}$ 
since the evaluations $\pmd{X}\ni\phi\mapsto\phi(x)$ ($\norm{x}\leq1$)
are completely contractive on $Z$ by Theorem \ref{acb}. This proves (i).

(ii) Given
$\phi=[\phi_{ij}]\in\matn{\pmd{(\minab{X})}}=
{\rm CB}_A(\minab{X},\matn{\bkh})_B$, the norm of $\phi$ is 
\begin{equation}\label{652} \norm{\phi}=\sup\{\norm{[\phi_{ij}(x_{kl})]}:\ 
[x_{kl}]\in{\rm M}_s(X),\ \m{[x_{kl}]}\leq1,\ s\in\bn\}.
\end{equation}
Since $\pmd{(\minab{X})}$ is a normal dual operator 
$\com{A},\com{B}$-bimodule, $\norm{\phi}\leq\Mnd{\phi}$
by maximality of $\Mnd{\cdot}$. For the reverse inequality, it suffices to show that
\begin{equation}\label{653} \norm{[T\phi_{ij}]}\leq\norm{\phi}
\end{equation}
for each contraction $T\in{\rm N}_{\com{A}}(\pmd{X},\matm{\bkh})_{\com{B}}$ 
($m\in\bn$). Let $\trm$ be the predual of $\matm{\bc}$ and put 
$Y=A\ehg\trm\ehg B$.
Since for each $n$ the unit ball of $\matn{A\hg\trm\hg B}$ is dense in that of
$\matn{Y}$ in the $A,B$-topology (by a similar argument as that 
preceding (\ref{444})), 
we have that
$\pmd{Y}=\acbbb{A\hg\trm\hg B}{\bkh}=\cbb{\trm}{\bkh}$, hence
\begin{equation}\label{654} \pmd{Y}=\pmd{(A\ehg\trm\ehg B)}=\matm{\bkh}.
\end{equation}

Realizing $X$ isometrically as a normal $A,B$-subbimodule in some ${\rm B}(\l)$,
let $\tilde{X}$ be the smallest strong $A,B$-bimodule  containing
$X$. Note that $\pmp{(\pmd{X})}=\tilde{X}$ by Theorem \ref{mpmd}, hence
$\pmd{\tilde{X}}=\pmd{X}$ by Theorem \ref{t38}. Since 
$\minab{X}\subseteq\minab{\tilde{X}}$ by (\ref{22}), it follows that 
replacing $X$ by $\tilde{X}$ has no effect on the statement (ii). In other words, we
may assume that $X$ is strong. We may regard $T$ as a complete contraction
from ${\rm MIN}_{\com{A}}(\pmd{X})_{\com{B}}$ into 
${\rm MIN}_{\com{A}}(\pmd{Y})_{\com{B}}$ (using (\ref{22}) for 
norms in $\matn{\pmd{Y}}$). Since these are normal dual operator 
bimodules by part (i), we deduce by Propositions \ref{p310} and \ref{p311}(i)
that $T=\pmd{S}$ for a contraction $S\in\B_A(Y,X)_B$.
Then the
norm of $[T\phi_{ij}]\in{\rm M}_{mn}(\bkh)
=\acbbb{A\ehg{\rm T}_{mn}\ehg B}{\bkh}$ (we have used (\ref{654}))
is equal to
\begin{equation}\label{655} \norm{[T\phi_{ij}]}=\sup\norm{[(T\phi_{ij})(v_{kl})]}
=\sup\norm{[\phi_{ij}(Sv_{kl})]},
\end{equation}
where the supremum is over all $[v_{kl}]\in{\rm M}_r(A\ehg {\rm T}_{mn}\ehg B)$
with $\norm{[v_{kl}]}\leq1$ and $r\in\bn$. Since $S$ is a 
complete contraction into
$\minab{X}$,  $\m{[Sv_{kl}]}\leq1$. Thus the right side
of (\ref{655}) is dominated by $\norm{\phi}$ (by (\ref{652})), which proves (\ref{653}).
\end{proof}

\begin{corollary}\label{c77} $\maxnd{X}$ is a normal dual operator
$A,B$-bimodule for each $X\in\andrb$.
\end{corollary}

\begin{proof} Let $X=\du{V}$. If $(x_{\nu})$ is a net in  the unit ball of 
$\matn{\minab{X}}$ converging to $x\in\matn{X}$ in the topology induced by
$\matn{V}$, then for each $a\in\rown{A}$ and $b\in\coln{B}$ the net $(ax_{\nu}b)$
converges to $axb$ in the topology induced by $V$. Since $X\in\andrb$, it follows
that $\norm{axb}\leq\norm{a}\norm{b}$ and, using (\ref{22}), we see that the
unit ball of $\matn{\minab{X}}$ is closed for each $n$.
By \cite[3.1]{Le} this implies
that $\minab{X}$ is a dual operator space and it follows that 
$\minab{X}\in\andob$ (using Theorem \ref{thbez} and Remark \ref{rnd}). Then by Theorem \ref{t38} (applied to
$\minab{X}$) we have in particular that 
$X=\pmd{(\pmp{X})}$ isometrically and weak* homeomorphically. Now 
Proposition \ref{p65}(ii) applied to $\pmp{X}$ shows that 
$\maxnd{X}=\maxnd{\pmd{(\pmp{X})}}=
\pmd{({\rm MIN}_{\com{A}}(\pmp{X})_{\com{B}})}$, which is a normal dual operator 
bimodule by Proposition \ref{pn}.
\end{proof}

\begin{example}\label{e67} In general $\minab{X^{\natural_{\rm p}\natural_{\rm p}}}\ne
(\minab{X})^{\natural_{\rm p}\natural_{\rm p}}$, hence by Proposition
\ref{p65} we have that
${\rm MAXND}_{\com{A}}(\pmd{X})_{\com{B}}=\pmd{(\minab{X})}\ne{\rm MAXN}_{\com{A}}(\pmd{X})_{\com{B}}.$

We sketch a counterexample. Let $A$ be the injective ${\rm II}_1$ factor represented normally on a Hilbert
space $\l$ such that $\l$ is not locally cyclic for $A$. Let 
$X=A\spac\com{A}\subseteq{\rm B}(\l\otimes\l)$. By \cite[3.4]{PSS}
$A$ (identified with $A\otimes 1$) is a norming subalgebra of $X$, which by (\ref{22}) means 
that $X$ carries the minimal operator $A$-bimodule structure.
By Corollary \ref{c35} $X^{\natural_{\rm p}\natural_{\rm p}}=\du{(X^{_A\sharp_A})}
\subseteq\bdu{X}=\tilde{X}$. Let $\g$ be the Hilbert space of the universal representation
$\Phi$ of $X$ (hence $\tilde{X}$ is the weak* closure of $\Phi(X)$). 
Since $A$ is a C$^*$-subalgebra of $X$, 
$\ta=A^{\sharp\sharp}$ can be regarded as a von Neumann subalgebra of
$\tilde{X}$. Let $P$ be the central projection in $\ta$ such that the weak*
continuous extension $\alpha$ of $\Phi^{-1}|\Phi(A)$ to $\ta$ has kernel 
$\ort{P}\ta$, so that $\alpha$ maps $P\ta$ isomorphically onto $A$.
Since $A$ is a factor, ${\rm C}^*(A\cup\com{A})$ is weak*
dense in $\bl$, hence the representation $a\otimes\com{a}
\mapsto a\com{a}$ of $X$ (bounded since $A$ is injective
\cite{EL}) is cyclic, therefore it can be regarded as a direct summand in
$\Phi$. So, we may regard $\l$ as a subspace in $\g$ and denote by $e\in\com{\tilde{X}}$
the projection onto $\l$. Then $\Phi(X)e\cong C^*(A\cup\com{A})$. If $C_e$ is
the central carrier of $e$ in $\tilde{X}$, the map
$$\tilde{X}C_e\to\tilde{X}e,\ \ x\mapsto xe$$
is an isomorphism of von Neumann algebras \cite[p. 335]{KR}, hence normal, and
maps the
C$^*$-subalgebra $\Phi(A\otimes1)C_e$ of $\tilde{X}C_e$ onto $\Phi(A\otimes1)e\cong A$.
Since the representation $a\mapsto\Phi(a\otimes1)|e\g$ of $A$ is just 
the identity, it is normal, hence the representation
$A\ni a\mapsto\Phi(a\otimes1)|C_e\g$ is also normal. Using \cite[10.1.18]{KR}
this implies that
$C_e\leq P$, hence $\tilde{X}C_e\subseteq P\tilde{X}P=P\bdu{X}P=
X^{\natural_{\rm p}\natural_{\rm p}}$ by Theorem \ref{t59}.  

If the operator $A$-bimodule structure on $X^{\natural_{\rm p}\natural_{\rm p}}$
were minimal, the same would hold for the subbimodule $\tilde{X}C_e$, hence also for
the completely isometric $A$-bimodule $\tilde{X}e$. But $\tilde{X}e\cong\bl$,
thus $\bl$ carries the
minimal operator $A$-bimodule structure, hence by (\ref{22}) $A$ is a norming
subalgebra of $\bl$. But this is a contradiction since by \cite[2.7]{PSS}
$A$ is norming for $\bl$ only if $\l$ is locally cyclic for $A$.
\end{example}

\begin{remark}\label{r68} By Proposition \ref{p65}(i)  and Corollary \ref{c35}(i)
${\rm MIN}_{\com{A}}(\pmd{X})_{\com{B}}$ is a dual operator space, hence
$\pmd{X}$ is the dual of a Banach space $V$. 
If there is an operator space  on $V$ such that $Y:
={\rm MAXN}_{\com{A}}(\pmd{X})_{\com{B}}$
is the operator space dual of $V$, then 
$Y$ is a normal {\em dual}  operator
$\com{A},\com{B}$-bimodule (Theorem \ref{thbez} and
Remark \ref{rnd}), hence 
$Y={\rm MAXND}_{\com{A}}(\pmd{X})_{\com{B}}$ by maximality. But, with $X$ as in  
Example \ref{e67}, $Y\ne{\rm MAXND}_{\com{A}}(\pmd{X})_{\com{B}}$, hence
there is no operator space on $V$ predual to $Y$.
\end{remark}

\end{document}